\documentclass[12pt]{article}
\usepackage{graphicx}
\usepackage{amsmath,amsthm,amssymb,enumerate}
\usepackage{euscript,mathrsfs}
\usepackage{color}
\usepackage{dsfont}
\usepackage[left=2cm,right=2cm,top=3.5cm,bottom=3.5cm]{geometry}
\usepackage{color}
\usepackage[framemethod=tikz]{mdframed}
\allowdisplaybreaks

\usepackage{soul}

\catcode`\@=11 \@addtoreset{equation}{section}

\catcode`\@=12

\allowdisplaybreaks

\newtheorem{Theorem}{Theorem}[section]
\newtheorem{Proposition}[Theorem]{Proposition}
\newtheorem{Lemma}[Theorem]{Lemma}
\newtheorem{Corollary}[Theorem]{Corollary}

\theoremstyle{definition}
\newtheorem{Definition}[Theorem]{Definition}

\newtheorem{Remark}[Theorem]{Remark}

\newcommand{\bTheorem}[1]{
\begin{Theorem} \label{T#1} }
\newcommand{\eT}{\end{Theorem}}

\newcommand{\bProposition}[1]{
\begin{Proposition} \label{P#1}}
\newcommand{\eP}{\end{Proposition}}

\newcommand{\bLemma}[1]{
\begin{Lemma} \label{L#1} }
\newcommand{\eL}{\end{Lemma}}

\newcommand{\bCorollary}[1]{
\begin{Corollary} \label{C#1} }
\newcommand{\eC}{\end{Corollary}}

\newcommand{\bRemark}[1]{
\begin{Remark} \label{R#1} }
\newcommand{\eR}{\end{Remark}}

\newcommand{\bDefinition}[1]{
\begin{Definition} \label{D#1} }
\newcommand{\eD}{\end{Definition}}

\newcommand{\vrE}{\vr_E}

\newcommand{\tmu}{\widetilde{\mu}}
\newcommand{\teta}{\widetilde{\eta}}
\newcommand{\tkap}{\widetilde{\kappa}}

\newcommand{\Ds}{\mathbb{D}_x}

\newcommand{\bFormula}[1]{
\begin{equation} \label{#1}}
\newcommand{\eF}{\end{equation}}

\newcommand{\Ov}[1]{\overline{#1}}

\newcommand{\aleq}{\stackrel{<}{\sim}}
\newcommand{\vd}{\vc{v}_\delta}

\newcommand{\vr}{\varrho}

\newcommand{\vt}{\vartheta}
\newcommand{\vu}{\vc{u}}
\newcommand{\vm}{\vc{m}}

\newcommand{\vq}{\vc{q}}

\newcommand{\vc}[1]{{\bf #1}}

\newcommand{\Div}{{\rm div}_x}
\newcommand{\Grad}{\nabla_x}

\newcommand{\dx}{\,{\rm d} {x}}

\newcommand{\dt}{\,{\rm d} t }

\newcommand{\vU}{\vc{U}}

\newcommand{\intO}[1]{\int_{\Omega} #1 \ \dx}

\newcommand{\D}{{\rm d}}

\newcommand{\ep}{\varepsilon}

\renewcommand{\S}{\mathbb{S}}

\newcommand{\br}{ \nonumber \\ }

\def\softd{{\leavevmode\setbox1=\hbox{d}%
          \hbox to 1.05\wd1{d\kern-0.4ex{\char039}\hss}}}
\definecolor{Cgrey}{rgb}{0.85,0.85,0.85}
\definecolor{Cblue}{rgb}{0.50,0.85,0.85}
\definecolor{Cred}{rgb}{1,0,0}
\definecolor{fancy}{rgb}{0.10,0.85,0.10}

\newcommand\Cbox[2]{%
    \newbox\contentbox%
    \newbox\bkgdbox%
    \setbox\contentbox\hbox to \hsize{%
        \vtop{
            \kern\columnsep
            \hbox to \hsize{%
                \kern\columnsep%
                \advance\hsize by -2\columnsep%
                \setlength{\textwidth}{\hsize}%
                \vbox{
                    \parskip=\baselineskip
                    \parindent=0bp
                    #2
                }%
                \kern\columnsep%
            }%
            \kern\columnsep%
        }%
    }%
    \setbox\bkgdbox\vbox{
        \color{#1}
        \hrule width  \wd\contentbox %
               height \ht\contentbox %
               depth  \dp\contentbox
        \color{black}
    }%
    \wd\bkgdbox=0bp%
    \vbox{\hbox to \hsize{\box\bkgdbox\box\contentbox}}%
    \vskip\baselineskip%
}

\mdfdefinestyle{MyFrame}{%
	linecolor=black,
	outerlinewidth=1pt,
	roundcorner=5pt,
	innertopmargin=\baselineskip,
	innerbottommargin=\baselineskip,
	innerrightmargin=10pt,
	innerleftmargin=10pt,
	backgroundcolor=gray!20!white}

\date{}

\begin{document}


\title{Euler system with a polytropic equation of state as a vanishing viscosity limit}

\author{Eduard Feireisl\thanks{The research of E.F. leading to these results has received funding from the
		Czech Sciences Foundation (GA\v CR), Grant Agreement
		21-02411S. The Institute of Mathematics of the Academy of Sciences of
		the Czech Republic is supported by RVO:67985840.} \and Christian Klingenberg\and Simon Markfelder\thanks{S. M. acknowledges financial support by the Alexander von Humboldt Foundation.}}

\date{\today}

\maketitle

\bigskip

\centerline{Institute of Mathematics of the Academy of Sciences of the Czech Republic}

\centerline{\v Zitn\' a 25, CZ-115 67 Praha 1, Czech Republic}

\medskip

\centerline{Institute of Mathematics, W\"urzburg University}

\centerline{Emil-Fischer-Str. 40, 97074 W\"urzburg, Germany}

\medskip

\centerline{Department of Applied Mathematics and Theoretical Physics, University of Cambridge} 

\centerline{Wilberforce Road, Cambridge CB3 0WA, United Kingdom}

\bigskip

\begin{abstract}
	
We consider the Euler system of gas dynamics endowed with the incomplete $(e - \vr - p)$ equation of state relating the internal energy $e$ to the mass density $\vr$ and the pressure $p$. 
We show that any sufficiently smooth solution can be recovered as a vanishing viscosity -- heat conductivity limit of the Navier--Stokes--Fourier system with a properly defined temperature. 
The result is unconditional in the case of the Navier type (slip) boundary conditions and extends to the no-slip condition for the velocity under some extra hypotheses of Kato's type 
concerning the behavior of the fluid in the boundary layer.

\end{abstract}


\bigskip


\section{Introduction}
\label{i}

The \emph{Euler system} describing the evolution of the density $\vr = \vr(t,x)$, the velocity $\vu = \vu(t,x)$, and the internal energy $e = e(t,x)$ of a compressible inviscid fluid reads
\begin{align} 
	\partial_t \vr + \Div (\vr \vu) &= 0, \br 
	\partial_t (\vr \vu) + \Div (\vr \vu \otimes \vu) + \Grad p &= 0,\br
	\partial_t \left[ \frac{1}{2} \vr |\vu|^2 + \vr e \right] + \Div\left(  \left( \left[ \frac{1}{2} \vr |\vu|^2 + \vr e \right] + p \right) \vu \right) &= 0.
	\label{i1}
	\end{align} 
The fluid is confined to a bounded domain $\Omega \subset R^3$, with \emph{impermeable boundary}, 
\begin{equation} \label{i2}  
\vu \cdot \vc{n} |_{\partial \Omega} = 0.
\end{equation}
The system \eqref{i1} rewritten in terms of the phase variables $(\vr, \vu, e)$ is symmetric hyperbolic, see e.g. Benzoni-Gavage and Serre \cite[Chapter 13, Section 13.2.2]{BenSer}.
The problem is formally closed by prescribing a suitable \emph{equation of state} (EOS). We consider a polytropic EOS 
\begin{equation} \label{i3}
	p = (\gamma - 1) \vr e \ \mbox{with the adiabatic exponent}\ \gamma > 1.
	\end{equation}

The equation of state \eqref{i3} is \emph{incomplete}, in particular, the (absolute) temperature $\vt$ is not uniquely determined. Indeed Gibbs' law asserts 
\begin{equation} \label{i4} 
	\vt D s = D e + p D \left( \frac{1}{\vr} \right), 
	\end{equation} 
where $s$ is a new thermodynamic variable called \emph{entropy}. Here $D=\left(\frac{\partial}{\partial \vr},\frac{\partial}{\partial \vt}\right)$. Plugging \eqref{i3} in \eqref{i4} we obtain a first order system that can be integrated yielding 
\begin{equation} \label{i5} 
	p(\vr, \vt) = \vt^{\frac{\gamma}{\gamma - 1}} P \left( \frac{\vr}{ \vt^{\frac{1}{\gamma - 1}}} \right), 
\end{equation}	
and, in accordance with \eqref{i3}, \eqref{i4}, 
\begin{align} 
e(\vr, \vt) &= \frac{\vt}{\gamma - 1}\frac{\vt^{\frac{1}{\gamma - 1}}}{\vr} P \left( \frac{\vr}{ \vt^{\frac{1}{\gamma - 1}}} \right), \br 
s(\vr, \vt) &= S \left( \frac{\vr}{ \vt^{\frac{1}{\gamma - 1}}} \right),\ S'(Z) = - \frac{1}{\gamma - 1} \frac{\gamma P(Z) - P'(Z) Z }{Z^2}, 
\label{i6}
\end{align}
for an arbitrary function $P$. Thus the absolute temperature $\vt$ is determined by $\vr$ and $e$ modulo the function $P$, see Cowperthwaite \cite{Cowp}, M\" uller and Ruggeri \cite{MURU}, or \cite[Chapters 2,3]{FeNo6A}.

The \emph{Navier--Stokes--Fourier system} describing the motion of a real viscous and heat conductive gas can be viewed as a viscous regularization of \eqref{i1}:
\begin{align} 
	\partial_t \vr + \Div (\vr \vu) &= 0, \br 
	\partial_t (\vr \vu) + \Div (\vr \vu \otimes \vu) + \Grad p &= \Div \mathbb{S},\br
	\partial_t (\vr s) + \Div (\vr s \vu) + \Div \left( \frac{ \vc{q} }{\vt} \right) &= \frac{1}{\vt} \left( \mathbb{S} : \Ds \vu - \frac{\vc{q} \cdot \Grad \vt }{\vt} \right) ,\ 
	\Ds \vu \equiv \frac{ \Grad \vu + \Grad^t \vu }{2},
	\label{i7}
\end{align} 
with the viscous stress $\mathbb{S}$ given by \emph{Newton's rheological law} 
\begin{equation} \label{i8}
	\mathbb{S} = 2 \tmu \left( \Ds \vu - \frac{1}{d} \Div \vu \mathbb{I} \right) + \teta \Div \vu \mathbb{I},
	\end{equation}
and the heat flux given by \emph{Fourier's law}
\begin{equation} \label{i9} 
	\vc{q} = - \tkap \Grad \vt.
	\end{equation} 
The Second law of thermodynamics requires the entropy production rate 
\[
\frac{1}{\vt} \left( \mathbb{S} : \Ds \vu - \frac{\vc{q} \cdot \Grad \vt }{\vt} \right)
\]
to be non--negative; whence the diffusion \emph{transport coefficients} $\tmu$, $\teta$, and $\tkap$ must be non--negative.
Note that, unlike in the Euler system \eqref{i1}, the knowledge of the temperature $\vt$ is necessary to determine the entropy as well as the heat flux in \eqref{i7}. The internal energy $e$ can be evaluated in terms of $\vr$, $\vt$ 
through \eqref{i6}. Thus solutions of the associated Navier--Stokes--Fourier system \eqref{i7}, that may be seen as a viscous regularization of the Euler system \eqref{i1}, depend on the choice of $P$ in \eqref{i5}.

We consider the vanishing dissipation limit of the Navier--Stokes--Fourier system, specifically, we rescale
\begin{equation} \label{i10}
	\mathbb{S}_n \approx \mu_n \mathbb{S},\ \vc{q}_n \approx \kappa_n \vc{q},
\ \mu_n \searrow 0,\ \kappa_n \searrow 0.
	\end{equation}
Moreover, 
the existing mathematical theory of the Navier--Stokes--Fourier system (see \cite{FeNo6A}) is based on the augmentation of the pressure, and, accordingly, the internal energy and entropy, by the radiation component 
\begin{equation} \label{i11}
	p_R = \frac{a}{3} \vt^4,\ e_R = \frac{a}{\vr} \vt^4,\ s_R = \frac{4a}{3 \vr} \vt^3,\ a > 0.
	\end{equation}
The parameter $a$ is very small and usually neglected in the real world applications. Consistently with \eqref{i10}, we therefore consider 
\begin{equation} \label{i12}
	a = a_n  ,\ a_n \searrow 0.
	\end{equation}

Suppose that $\gamma > 1$ is given and that the Euler system \eqref{i1}--\eqref{i3} admits a smooth ($C^1$) solution on a time interval $[0, T]$. Our goal is to identify the function $P$ in \eqref{i5} in such a way that any sequence of 
weak solutions to the Navier--Stokes--Fourier system \eqref{i7}--\eqref{i9} converges in the vanishing viscosity/radiation limit \eqref{i10}--\eqref{i12} to the solution of the Euler system in $(0,T) \times \Omega$. 
Moreover, we show that
the convergence 
is unconditional, if the boundary layer is eliminated by the choice of the complete slip boundary conditions 
\begin{equation} \label{i13}
	\vu \cdot \vc{n}|_{\partial \Omega} = 0,\ (\mathbb{S} \cdot \vc{n}) \times \vc{n}|_{\partial \Omega} = 0,
	\end{equation}
where $\vc{n}$ denotes the outer normal vector to $\partial \Omega$. In the case of the no--slip boundary conditions 
\begin{equation} \label{i14}
	\vu|_{\partial \Omega} = 0, 
	\end{equation}
the convergence is conditioned by extra hypotheses of Kato's type \cite{Kato}, \cite{Kato1984} identified in the compressible setting by Sueur \cite{Sue1} and Wang, Zhu \cite{WangZhu}.

In comparison with the existing literature, notably \cite{WangZhu}, our result covers all admissible values of the adiabatic coefficient $\gamma$ in \eqref{i3} as well as general dependence of the transport coefficients 
on the temperature in the spirit of the existence theory developed in \cite{FeNo6A}. 

The paper is organized as follows. In Section \ref{p}, we recall the necessary preliminary material concerning the weak solutions to the Navier--Stokes--Fourier system including the relative energy inequality that represents a crucial tool in the analysis. Section \ref{r} contains the main results. In Section \ref{c} we show consistency of the vanishing viscosity approximation. Specifically, the viscous stress, the heat flux as well as the radiation components of the 
pressure, internal energy, and entropy along with the associated fluxes disappear in the regime specified in \eqref{i10}, \eqref{i12}. This process is ``path dependent'', specifically certain relations concerning the asymptotic behaviour of $(\mu_n,  \kappa_n, a_n)$ must be imposed in the spirit of \cite{Fei2015A}. The convergence towards the strong solution of the Euler system is shown in Section \ref{C}.

\section{Preliminary material}
\label{p} 

We recall the existing theory of weak solutions to the Navier--Stokes--Fourier system. 

\subsection{Mathematical theory of the closed system}

By imposing either the complete slip \eqref{i13} or the no slip boundary condition \eqref{i14} we suppose the fluid is mechanically insulated. In view of our final objective, we require the fluid to be energetically 
isolated, specifically 
\begin{equation} \label{p1}
	\vc{q} \cdot \vc{n}|_{\partial \Omega} = 0.
	\end{equation}
The mathematical theory for closed systems relevant for future analysis was developed in \cite{FeNo6A}. Note that the extension to \emph{open} systems is also available in the recent works \cite{ChauFei} and \cite{FeiNov20}, see also the forthcoming monograph \cite{FeiNovOpen}. 

A suitable \emph{weak formulation} of the Navier--Stokes--Fourier system augmented by the radiative terms proposed in \cite{FeNo6A} reads
\begin{align} 
\partial_t \vr + \Div (\vr \vu) &= 0, \br 
\partial_t (\vr \vu) + \Div (\vr \vu \otimes \vu) + \Grad (p + p_R) &= \mu \Div \mathbb{S},\br	
\partial_t (\vr (s + s_R)) + \Div( \vr (s + s_R) \vu) + \kappa \Grad \left( \frac{ \vc{q} }{\vt} \right) &\geq \frac{1}{\vt} \left( \mu \mathbb{S} : \Ds \vu - \kappa \frac{\vc{q} \cdot \Grad \vt}{\vt} \right), \br
\frac{\D }{\dt} \intO{ \left( \frac{1}{2} \vr |\vu|^2 + \vr (e + e_R) \right) } &= 0,
\label{p2}
\end{align}              
see \cite[Chapter 3]{FeNo6A}. Note that we anticipate the influence of thermal radiation represented by the extra terms 
$p_R$, $e_R$, and $s_R$ in \eqref{p2}. In accordance with \eqref{i12}, these terms will vanish in the asymptotic limit. The energy balance appearing in the Euler system in \eqref{i1} is replaced by 
the entropy inequality supplemented with the total energy balance in \eqref{p2}. 

\subsection{Transport coefficients}

In accordance with the molecular theory of gases (see e.g. Becker \cite{BE}), the transport coefficients depend on the temperature. Specifically, we assume that $\tmu$, $\teta$, $\tkap$ are continuously differentiable functions of $\vt$ satisfying
\begin{align} 
	0 < \underline{\mu} \left( 1  + \vt^\alpha \right) &\leq \tmu(\vt) \leq \Ov{\mu}  \left( 1  + \vt^\alpha \right),\ \alpha \geq 0, \br
	\sup_{ \vt \in [0, \infty) }|\tmu'(\vt)| &< \infty,\ \br
	0 \leq \teta (\vt) \leq  \Ov{\eta}  \left( 1  + \vt^\alpha \right),\
	0 < \underline{\kappa} \left( 1  + \vt^3 \right) &\leq \tkap(\vt) \leq \Ov{\kappa}  \left( 1  + \vt^3 \right)
	\label{p3}
	\end{align}
for any $\vt \geq 0$. Note that the cubic growth of $\kappa$ is motivated by the presence of the radiation terms, see Oxenius \cite{OX}. 

\subsection{Equation of state}

A proper choice of the equation of state for the Navier--Stokes--Fourier system plays of course a crucial role in the present paper. Given $\gamma > 1$, we have to identify the function $P$ in \eqref{i5}. For $p = p(\vr, \vt)$, 
$e = e(\vr, \vt)$, we recall the \emph{hypothesis of thermodynamic stability} 
\begin{equation} \label{p4}
	\frac{\partial p(\vr, \vt) }{\partial \vr} > 0,\ 	\frac{\partial e(\vr, \vt) }{\partial \vt} > 0.
\end{equation}
This imposes the following restrictions on $P$: 
\begin{align} 
	P'(Z) &> 0 \ \mbox{for all}\ Z > 0 , \br
	\gamma P(Z) - P'(Z) Z  &> 0 \ \mbox{for all}\ Z > 0.
	\label{p5}
\end{align}

The following lemma shows existence of a suitable $P$. 

\begin{Lemma} \label{lemma:P}
	For all $\underline{Z}>0$ there exist functions $P,S\in C^1[0,\infty)$ with properties \eqref{i6}, \eqref{p5} and such that 
	\begin{equation} \label{eq:lemmaP}
		P(Z) = Z  \quad \text{ for all }Z\in [0,\underline{Z}].
	\end{equation}
	Moreover $P,S$ satisfy
	\begin{align}
		P(0) &= 0 ; \label{eq:prop-SP1}\\
		\frac{\gamma P(Z) - P'(Z) Z}{Z} &\leq C \quad \text{ for all } Z> 0 ; \label{eq:prop-SP2} \\ 
		\lim_{Z\to\infty} \frac{P(Z)}{Z^\gamma} &>0; \label{eq:prop-SP3}  \\
		\lim\limits_{Z\to \infty} S(Z) &=0 .  \label{eq:prop-SP4}
	\end{align} 
\end{Lemma}

Note that according to \eqref{eq:prop-SP4}, $S$ from lemma \ref{lemma:P} is in accordance with the Third law of thermodynamics, namely 
\[
s(\vr, \vt) \to 0 \ \mbox{as}\ \vt \to 0+ \ \mbox{for any fixed}\ \vr > 0, 
\] 
cf. Belgiorno \cite{BEL1}, \cite{BEL2}.

\begin{proof}
	Let us first consider the case $\underline{Z}=1$. Set $P,S\in C^1[0,\infty)$ 
	$$
	P(Z):= \left\{ \begin{array}{ll}
		Z & \text{ if }Z\leq 1 ,\\
		\frac{\gamma-1}{\gamma} + \frac{1}{\gamma} Z^\gamma & \text{ if } Z>1 ,
	\end{array} \right.
	$$
	and 
	$$
	S(Z):= \left\{ \begin{array}{ll}
		-\log(Z) + 1 & \text{ if }Z\leq 1 ,\\
		\frac{1}{Z} & \text{ if } Z>1 ,
	\end{array} \right.
	$$
	It is then straightforward to check \eqref{i6}, \eqref{p5}, \eqref{eq:lemmaP}--\eqref{eq:prop-SP4}. 
	
	Let us now look at $\underline{Z}\neq 1$. We define the $P,S$ constructed above as $P_1,S_1$ and set 
	$$
	P(Z):= P_1 \left(\frac{Z}{\underline{Z}}\right), \qquad S(Z):= \frac{1}{\underline{Z}} S_1 \left(\frac{Z}{\underline{Z}}\right).
	$$
	Again straightforward computations show that the properties \eqref{i6}, \eqref{p5}, \eqref{eq:lemmaP}--\eqref{eq:prop-SP4} follow from the corresponding property of $P_1,S_1$. 
\end{proof} 

Note that for $Z\in [0,\underline{Z}]$, according to \eqref{eq:lemmaP} and \eqref{i5} we simply obtain the Boyle-Mariotte law 
$$
p(\vr, \vt) = \vr \vt .
$$
Hence the temperature for the Euler system \eqref{i1} endowed with the incomplete EOS \eqref{i3} can be recovered by choosing $\underline{Z}$ in lemma \ref{lemma:P} appropriately, see section \ref{subsec:temp} for details. 

\subsection{Relative energy}

The \emph{relative energy} for the Navier--Stokes--Fourier system may be seen as a counterpart of Dafermos' relative entropy for the (hyperbolic) Euler system, see \cite{Daf4}. Given a trio of ``test functions''
\begin{equation} \label{p9}
r > 0, \ \Theta > 0,\ \vc{U} ,
\end{equation}
the relative energy reads 
\begin{align} 
	E \left( \vr, \vt, \vu \Big| r, \Theta, \vc{U} \right) = 
	\frac{1}{2} \vr |\vu - \vc{U} |^2 + H_\Theta (\vr, \vt) - \frac{\partial H_\Theta (r,\Theta) }{\partial \vr} (\vr - r) - H_\Theta (r, \Theta)
	\label{p10}
	\end{align}
where 
\[
H_\Theta (\vr, \vt) = \vr (e(\vr, \vt) - \Theta s(\vr, \vt))
\]
is the ballistic free energy. In the context of the system \eqref{p2} perturbed by the radiation terms, we have 
\[
H_\Theta (\vr, \vt) = \vr ((e + e_R)(\vr, \vt) - \Theta (s + s_R) (\vr, \vt)).
\]
The relative energy augmented by the radiation component will be denoted $E_a$.
We also introduce the standard energy
\[
E(\vr, \vt, \vu) = \frac{1}{2} \vr |\vu|^2 + \vr e(\vr, \vt).
\]

The following result was proved in \cite{FeiNov10}: 
Suppose that:

\begin{itemize}
\item $(\vr, \vt, \vu)$ is a weak solution to the Navier--Stokes--Fourier system \eqref{p2} in $(0,T) \times \Omega$ with the no--flux boundary conditions \eqref{p1} and either the complete slip boundary conditions 
\eqref{i13} or the no--slip boundary condition \eqref{i14}.

\item $(r, \Theta, \vc{U})$ is a trio of continuously differentiable test functions, 
\[
r > 0,\ \Theta > 0 \ \mbox{in}\ [0,T] \times \Ov{\Omega}, 
\]
where $\vc{U}$ satisfies either the impermeability boundary condition 
\[
\vc{U} \cdot \vc{n}|_{\partial \Omega} = 0, 
\]
or the no--slip boundary condition 
\[
\vc{U}|_{\partial \Omega} = 0.
\]
	
	\end{itemize}

Then the \emph{relative energy inequality}

\begin{align}
	&\Big[ \intO{ E_a \Big(\vr,\vt,\vu\big| r, \Theta, \vU \Big) } \Big]_{t=0}^{t=\tau} + \int_0^\tau \int_\Omega \frac{\Theta}{\vt}\left( \mu \S(\vt,\Grad \vu):\Grad\vu - \kappa \frac{\vq(\vt,\Grad\vt)\cdot \Grad\vt}{\vt}\right)\dx\dt \notag\\
	&\leq\int_0^\tau \int_\Omega \vr(\vu-\vU)\cdot\Grad \vU\cdot (\vU-\vu)\dx\dt \notag \\ &\quad + \mu \int_0^\tau \int_\Omega  \S(\vt,\Grad\vu) :\Grad\vU \dx\dt 
	- \kappa \int_0^\tau \int_\Omega \frac{\vq(\vt,\Grad\vt)}{\vt}\cdot \Grad\Theta \dx\dt \notag \\ &\quad+ \int_0^\tau \int_\Omega \vr \big((s + s_R) (\vr,\vt) - (s + s_R)(r,\Theta)\big)(\vU-\vu)\cdot \Grad\Theta \dx\dt \notag \\
	&\quad +\int_0^\tau\int_\Omega \vr \Big(\partial_t\vU + \vU\cdot \Grad \vU\Big)\cdot(\vU-\vu)\dx\dt - \int_0^\tau\int_\Omega (p + p_R)(\vr,\vt)\Div\vU \dx\dt \notag\\
	&\quad - \int_0^\tau\int_\Omega \vr \big((s + s_R) (\vr,\vt)-(s + s_R) (r,\Theta)\big)\Big(\partial_t\Theta + \vU\cdot \Grad\Theta\Big)\dx\dt \notag\\ 
	&\quad + \int_0^\tau\int_\Omega \left(\left(1-\frac{\vr}{r}\right)\partial_t (p + p_R) (r,\Theta) - \frac{\vr}{r}\vu\cdot \Grad (p + p_R) (r,\Theta)\right)\dx\dt \label{p11}
\end{align} 
holds for a.a. $\tau \in (0,T)$.

Finally, we recall the fundamental properties of the relative energy that follow from the hypothesis of thermodynamic stability \eqref{p4}. In accordance with hypothesis \eqref{p9}, fix 
\begin{align}
0 < \underline{\vr} &< \inf_{[0,T] \times \Ov{\Omega}} r \leq \sup_{[0,T] \times \Ov{\Omega}} r < \Ov{\vr}, \br
0 < \underline{\vt} &< \inf_{[0,T] \times \Ov{\Omega}} \Theta \leq \sup_{[0,T] \times \Ov{\Omega}} \Theta < \Ov{\vt}, 
\nonumber
\end{align}
and define 
\[
[F]_{\rm ess} =  \Phi (\vr, \vt) F,\ 
[F]_{\rm res} = F - [F]_{\rm ess}, 
\]
where 
\[
\Phi \in C^1_c(0,\infty)^2,\ 0 \leq \Phi \leq 1,\ \Phi (\vr, \vt) = 1 \ \mbox{whenever}\ \underline{\vr}  \leq \vr \leq \Ov{\vr} \ \mbox{and}\ \underline{\vt} \leq \vt \leq\Ov{\vt}.
\]

Then
\begin{align}
	E_a \left(\vr, \vt, \vu \Big| r, \Theta, \vu \right) \geq 
	E \left(\vr, \vt, \vu \Big| r, \Theta, \vu \right) &\geq c \left( [\vr- r ]_{\rm ess}^2  + [\vt - \Theta]_{\rm ess}^2 + [\vu - \vU ]^2_{\rm ess}  \right) \br
	E_a \left(\vr, \vt, \vu \Big| r, \Theta, \vu \right) &\geq c \left( 1_{\rm res} + [\vr (e + e_R)(\vr, \vt)]_{\rm res} + [ \vr |(s + s_R)(\vr, \vt)| ]_{\rm res} \right) ,\br 
	E \left(\vr, \vt, \vu \Big| r, \Theta, \vu \right) &\geq c \left( 1_{\rm res} + [\vr e (\vr, \vt)]_{\rm res} + [ \vr | s (\vr, \vt)| ]_{\rm res} \right),
		\label{p12}
	\end{align}
where the constants depend on $\underline{\vr}$, $\Ov{\vr}$, $\underline{\vt}$, and $\Ov{\vt}$, see e.g. \cite{FeNo6A} for details. As a consequence of the hypothesis of thermodynamic stability \eqref{p4}, the relative energy expressed in terms of the 
conservative entropy variables $(\vr, \vc{m} = \vr \vu, \mathcal{S} = \vr s )$ is a strictly convex function and represents to so-called Bregman distance between 
$(\vr, \vm, \mathcal{S})$ and $(r, r\vU, r s(r, \Theta) )$, see e.g. \cite{FeiNov20}.
Note carefully that the relative energy $E_a$ associated to the Navier--Stokes--Fourier system \eqref{p2} is augmented by the radiation component
\[
a (\vt^4 - \Theta^4) + \frac{4a}{3} \Theta (\Theta^3 - \vt^3 ) \geq 0.
\]

\section{Main results}
\label{r}

We state the main results in the physically relevant case $\Omega \subset R^3$. 
We consider three vanishing parameters in the asymptotic limit: the viscosity coefficient $\mu_n$, the heat conductivity coefficient $\kappa_n$, and the 
radiation parameter $a_n$, cf. \eqref{i10}, \eqref{i12}.

\subsection{Unconditional convergence in the absence of boundary layer}

We start with the Navier--Stokes--Fourier system \eqref{p2}, with the complete slip boundary conditions \eqref{i13}, and the no--flux boundary condition \eqref{p1}.

\begin{Theorem} [Unconditional convergence] \label{Tr1}
	Let $\Omega \subset R^3$ be a bounded Lipschitz domain. Suppose that the Euler system \eqref{i1}--\eqref{i3}, with $\gamma > 1$, admits a strong solution 
	\[
	\vr_E, \ e_E \in C^1([0,T] \times \Ov{\Omega}), \ \vu_E \in C^1([0,T] \times \Omega; R^3)
	\]
	satisfying
	\[
	\inf_{[0,T] \times \Ov{\Omega} } \vr_E > 0,\ \inf_{[0,T] \times \Ov{\Omega} } e_E > 0.
	\]
	
	Then there exists a ($p$-$\vr$-$\vt$) EOS $p = p(\vr, \vt)$ that complies with Gibbs' relation \eqref{i4} as well as the hypothesis of thermodynamic stability \eqref{p4}, with the associated internal energy EOS $e = e(\vr, \vt)$
	and entropy $s = s(\vr, \vt)$ determined through \eqref{i6}, such that the following holds:
	
	Let $(\vr_n, \vt_n, \vu_n)_{n=1}^\infty$ be a sequence of weak solutions to the Navier--Stokes--Fourier system \eqref{p2}, with the complete slip boundary condition \eqref{i13}, and the no--flux boundary conditions 
	\eqref{p1}, in the vanishing dissipation/radiation regime:
	\begin{equation} \label{r1}
		\mu_n \searrow 0,\ a_n \approx \mu_n^{\frac{4}{1 + \alpha}}, \ \frac{\kappa_n}{a_n^{\frac{3}{4}}} \to 0,
		\end{equation}
	where $\alpha \in [\frac{1}{3}, 1]$ is the exponent in hypothesis \eqref{p3}.
	In addition, suppose that the initial data
	\[
	\vr_{n,0} = \vr_n(0, \cdot),\ \vt_{n,0} = \vt_n (0, \cdot),\ \vu_{n,0} = \vu_n(0, \cdot),
	\]
	converge strongly to those of the Euler system, specifically,
	\begin{align}
	0 &< \underline{\vr} < \inf_{(0,T) \times \Omega} \vr_{n,0} \leq  \sup_{(0,T) \times \Omega} \vr_{n,0} < \Ov{\vr}\ \mbox{uniformly in}\ n, \ \vr_{n,0} \to \vr_E(0, \cdot) \ \mbox{in}\ L^1 (\Omega), \br
	0 &< \underline{\vt} < \inf_{(0,T) \times \Omega} \vt_{n,0} \leq  \sup_{(0,T) \times \Omega} \vt_{n,0} < \Ov{\vt}\ \mbox{uniformly in}\ n,\ 
	e(\vr_{n,0}, \vt_{n,0})(0, \cdot) \to e_E(0, \cdot) \ \mbox{in}\ L^1(\Omega), \br
	&|\vu_{n,0} | \leq \Ov{\vu} \ \mbox{uniformly in}\ n,\ \vu_{n,0} \to \vu_E(0, \cdot) \ \mbox{in}\ L^1(\Omega; R^3).
\label{r2}
	\end{align}

Then 
\begin{equation} \label{r3}
\vr_n \to \vr_E,\ \vr_n e(\vr_n, \vt_n) \to \vr_E e_E  \ \mbox{in}\ L^1((0,T) \times \Omega),\ \vr_n \vu_n \to \vr_E \vu_E \ \mbox{in}\ L^1((0,T) \times \Omega; R^3). 
\end{equation}

		\end{Theorem}
	
	\begin{Remark} \label{Rr1} 
		The reader may consult \cite[Chapter 3]{FeNo6A} for the exact definition of a weak solution of the Navier--Stokes--Fourier system emanating from the initial data $(\vr_{n,0}, \vt_{n,0}, \vu_{n,0})$.
		\end{Remark} 
	
	\begin{Remark} \label{Rr2}
		We strongly point out that Theorem \ref{Tr1} \emph{does not} contain any claim concerning the \emph{existence} of weak solutions for the Navier--Stokes--Fourier system. The existence is known only in some particular cases: 
		$\gamma \geq \frac{5}{3}$, $\alpha \in [\frac{2}{5}, 1]$, see \cite[Chapter 3, Theorem 3.1]{FeNo6A}, and $\gamma > \frac{3}{2}$, $\alpha = 1$, see Jessl\' e, Jim, and Novotn\' y \cite[Theorem 2.1]{JesJiNov}.
		The best known results for the planar flows were obtained recently by Pokorn\'{y}, M. and Sk\v{r}\'{\i}\v{s}ovsk\'{y} \cite{PokSkr}.

		\end{Remark}
	Local in time existence of smooth solutions to the Euler system was established by Schochet \cite{SCHO1}.	
	
\subsection{Conditional result -- viscous boundary layer}	

The no--slip boundary condition \eqref{i14} imposed on the viscous flow cannot be retained for the limit Euler system and the well known problem of viscous boundary layer appears. We report a conditional result \` a la Kato in the spirit of Sueur \cite{Sue1} and Wang, Zhu \cite{WangZhu}. 
Let 
\[
\Omega_\delta = \left\{ x \in \Omega \ \Big|\  {\rm dist}[x, \partial \Omega]  < \delta \right\}.
\]
Any vector field $\vc{w}$ can be decomposed into its normal and tangential component with respect to $\partial \Omega$:
\begin{align}
	\vc{w}(t,x) &= \vc{w}_n(t,x) + \vc{w}_\tau(t,x), \br \vc{w}_n(t,x) &= \left( \vc{w} \cdot \Grad {\rm dist}[x, \partial \Omega] \right) \Grad {\rm dist}[x, \partial \Omega],\ 
	\vc{w}_\tau(t,x) = \vc{w}(t,x) - \vc{w}_n(t,x).
	\nonumber
\end{align}

Note that $|\Grad {\rm dist}[x, \partial \Omega]|=1$, see section \ref{subsection:C1}.

We start with a result inspired by Sueur \cite{Sue1}. 

\begin{Theorem} [Conditional convergence, gradient criterion] \label{Tr3}
	Let $\Omega \subset R^3$ be a bounded domain of class $C^{2 + \nu}$. Suppose that the Euler system \eqref{i1}--\eqref{i3}, with $\gamma > 1$, admits a strong solution 
	\[
	\vr_E, \ e_E \in C^1([0,T] \times \Ov{\Omega}), \ \vu_E \in C^1([0,T] \times \Omega; R^3)
	\]
	satisfying
	\[
	\inf_{[0,T] \times \Ov{\Omega} } \vr_E > 0,\ \inf_{[0,T] \times \Ov{\Omega} } e_E > 0.
	\]
	
	Then there exists a ($p$-$\vr$-$\vt$) EOS $p = p(\vr, \vt)$ that complies with Gibbs' relation \eqref{i4} as well as the hypothesis of thermodynamic stability \eqref{p4}, with the associated internal energy EOS $e = e(\vr, \vt)$
	and entropy $s = s(\vr, \vt)$ determined through \eqref{i6}, such that the following holds:
	
	Let $(\vr_n, \vt_n, \vu_n)_{n=1}^\infty$ be a sequence of weak solutions to the Navier--Stokes--Fourier system \eqref{p2}, with the no--slip boundary condition \eqref{i14}, and the no--flux boundary conditions 
	\eqref{p1}, in the vanishing dissipation/radiation regime:
	\[
		\mu_n \searrow 0,\ a_n \approx \mu_n^{\frac{4}{1 + \alpha}}, \ \frac{\kappa_n}{a_n^{\frac{3}{4}}} \to 0,
	\]
	where $\alpha \in [\frac{1}{3}, 1]$ is the exponent in hypothesis \eqref{p3}. In addition, suppose that the initial data
	\[
	\vr_{n,0} = \vr_n(0, \cdot),\ \vt_{n,0} = \vt_n (0, \cdot),\ \vu_{n,0} = \vu_n(0, \cdot),
	\]
	converge strongly to those of the Euler system, specifically,
	\begin{align}
		0 &< \underline{\vr} < \inf_{(0,T) \times \Omega} \vr_{n,0} \leq  \sup_{(0,T) \times \Omega} \vr_{n,0} < \Ov{\vr}\ \mbox{uniformly in}\ n, \ \vr_{n,0} \to \vr_E(0, \cdot) \ \mbox{in}\ L^1 (\Omega), \br
		0 &< \underline{\vt} < \inf_{(0,T) \times \Omega} \vt_{n,0} \leq  \sup_{(0,T) \times \Omega} \vt_{n,0} < \Ov{\vt}\ \mbox{uniformly in}\ n,\ 
		e(\vr_{n,0}, \vt_{n,0})(0, \cdot) \to e_E(0, \cdot) \ \mbox{in}\ L^1(\Omega), \br
		&|\vu_{n,0} | \leq \Ov{\vu} \ \mbox{uniformly in}\ n,\ \vu_{n,0} \to \vu_E(0, \cdot) \ \mbox{in}\ L^1(\Omega; R^3).
		\nonumber
	\end{align}
	Finally, suppose
	\begin{align}
		\mu_n &\int_0^T \int_{\Omega_{\mu_n}} | \mathbb{S} (\vt_n, \Grad \vu_n ) |^2 \dx \dt \to 0, \br 
		\mu_n &\int_0^T \int_{\Omega_{\mu_n}} \left( \frac{\vr_n |\vu_n|^2 }{{\rm dist}^2[x, \partial \Omega] } + \frac{\vr^2_n |(\vu_n)_n|^2 }{{\rm dist}^2[x, \partial \Omega] }
		\right) \dx \dt \to 0
		\label{rr2b}
	\end{align}
	as $n \to \infty$.	
	
	Then 
	\[
	\vr_n \to \vr_E,\ \vr_n e(\vr_n, \vt_n) \to \vr_E e_E  \ \mbox{in}\ L^1((0,T) \times \Omega),\ \vr_n \vu_n \to \vr_E \vu_E \ \mbox{in}\ L^1((0,T) \times \Omega; R^3). 
	\]
	
\end{Theorem}

Finally, we state a conditional result inspired by Wang and Zhu \cite{WangZhu}.

\begin{Theorem} [Conditional convergence] \label{Tr2}
	Let $\Omega \subset R^3$ be a bounded domain of class $C^{2 + \nu}$. Suppose that the Euler system \eqref{i1}--\eqref{i3}, with $\gamma > 1$, admits a strong solution 
	\[
	\vr_E, \ e_E \in C^1([0,T] \times \Ov{\Omega}), \ \vu_E \in C^1([0,T] \times \Omega; R^3)
	\]
	satisfying
	\[
	\inf_{[0,T] \times \Ov{\Omega} } \vr_E > 0,\ \inf_{[0,T] \times \Ov{\Omega} } e_E > 0.
	\]
	
	Then there exists a ($p$-$\vr$-$\vt$) EOS $p = p(\vr, \vt)$ that complies with Gibbs' relation \eqref{i4} as well as the hypothesis of thermodynamic stability \eqref{p4}, with the associated internal energy EOS $e = e(\vr, \vt)$
	and entropy $s = s(\vr, \vt)$ determined through \eqref{i6}, such that the following holds:
	
	Let $(\vr_n, \vt_n, \vu_n)_{n=1}^\infty$ be a sequence of weak solutions to the Navier--Stokes--Fourier system \eqref{p2}, with the no--slip boundary condition \eqref{i14}, and the no--flux boundary conditions 
	\eqref{p1}, in the vanishing dissipation/radiation regime:
	\[
		\mu_n \searrow 0,\ a_n \approx \mu_n^{\frac{4}{1 + \alpha}}, \ \frac{\kappa_n}{a_n^{\frac{3}{4}}} \to 0,
	\]
	where $\alpha \in [\frac{1}{3}, 1]$ is the exponent in hypothesis \eqref{p3}. In addition, suppose that the initial data
	\[
	\vr_{n,0} = \vr_n(0, \cdot),\ \vt_{n,0} = \vt_n (0, \cdot),\ \vu_{n,0} = \vu_n(0, \cdot),
	\]
	converge strongly to those of the Euler system, specifically,
	\begin{align}
		0 &< \underline{\vr} < \inf_{(0,T) \times \Omega} \vr_{n,0} \leq  \sup_{(0,T) \times \Omega} \vr_{n,0} < \Ov{\vr}\ \mbox{uniformly in}\ n, \ \vr_{n,0} \to \vr_E(0, \cdot) \ \mbox{in}\ L^1 (\Omega), \br
		0 &< \underline{\vt} < \inf_{(0,T) \times \Omega} \vt_{n,0} \leq  \sup_{(0,T) \times \Omega} \vt_{n,0} < \Ov{\vt}\ \mbox{uniformly in}\ n,\ 
		e(\vr_{n,0}, \vt_{n,0})(0, \cdot) \to e_E(0, \cdot) \ \mbox{in}\ L^1(\Omega), \br
		&|\vu_{n,0} | \leq \Ov{\vu} \ \mbox{uniformly in}\ n,\ \vu_{n,0} \to \vu_E(0, \cdot) \ \mbox{in}\ L^1(\Omega; R^3).
		\nonumber
	\end{align}
Finally, suppose there is a sequence $\delta_n \to 0$ such that 
\begin{align}
	\frac{\mu_n}{\delta_n} &\to 0 \ \mbox{as}\ n \to \infty, \br
	\frac{1}{\delta_n} &\int_0^T \int_{\Omega_{\delta_n}} \vt^{1 + \alpha}_n \dx \dt \leq c,\br 
	\int_0^T & \left( \frac{1}{\delta_n} \left\| \vr_n (\vu_n)_n \right\|_{L^{\frac{24}{17 + 3\alpha}}(\Omega_{\delta_n})} 
	+ \frac{1}{\delta^2_n \mu_n}	 \left\| \vr_n (\vu_n)_n \right\|_{L^{\frac{24}{17 + 3\alpha}}(\Omega_{\delta_n})}^2 \left\| \vt_n^{\frac{1 - \alpha}{2}} \right\|_{L^{\frac{8}{1 - \alpha}}(\Omega_{\delta_n})}^2 \right) \dt \to 0
	\label{r2b}
	\end{align}
uniformly for $n \to \infty$.	
	
	Then 
	\[
		\vr_n \to \vr_E,\ \vr_n e(\vr_n, \vt_n) \to \vr_E e_E  \ \mbox{in}\ L^1((0,T) \times \Omega),\ \vr_n \vu_n \to \vr_E \vu_E \ \mbox{in}\ L^1((0,T) \times \Omega; R^3). 
	\]
	
\end{Theorem}

Hypothesis \ref{r2b} is awkward and seems much stronger than its counter part by Wang and Zhu \cite{WangZhu}. However, for $\alpha = 1$ we have the following. 

\begin{Theorem} [Conditional convergence, $\alpha = 1$] \label{Tr2bis}
	Let $\Omega \subset R^3$ be a bounded domain of class $C^{2 + \nu}$. Suppose that the Euler system \eqref{i1}--\eqref{i3}, with $\gamma > 1$, admits a strong solution 
	\[
	\vr_E, \ e_E \in C^1([0,T] \times \Ov{\Omega}), \ \vu_E \in C^1([0,T] \times \Omega; R^3)
	\]
	satisfying
	\[
	\inf_{[0,T] \times \Ov{\Omega} } \vr_E > 0,\ \inf_{[0,T] \times \Ov{\Omega} } e_E > 0.
	\]
	
	Then there exists a ($p$-$\vr$-$\vt$) EOS $p = p(\vr, \vt)$ that complies with Gibbs' relation \eqref{i4} as well as the hypothesis of thermodynamic stability \eqref{p4}, with the associated internal energy EOS $e = e(\vr, \vt)$
	and entropy $s = s(\vr, \vt)$ determined through \eqref{i6}, such that the following holds:
	
	Let $(\vr_n, \vt_n, \vu_n)_{n=1}^\infty$ be a sequence of weak solutions to the Navier--Stokes--Fourier system \eqref{p2}, with the no--slip boundary condition \eqref{i14}, and the no--flux boundary conditions 
	\eqref{p1}, in the vanishing dissipation/radiation regime with $\alpha = 1$:
	\[
		\mu_n \searrow 0,\ a_n \approx \mu_n^2, \ \frac{\kappa_n}{a_n^{\frac{3}{4}}} \to 0.
	\]
	In addition, suppose that the initial data
	\[
	\vr_{n,0} = \vr_n(0, \cdot),\ \vt_{n,0} = \vt_n (0, \cdot),\ \vu_{n,0} = \vu_n(0, \cdot),
	\]
	converge strongly to those of the Euler system, specifically,
	\begin{align}
		0 &< \underline{\vr} < \inf_{(0,T) \times \Omega} \vr_{n,0} \leq  \sup_{(0,T) \times \Omega} \vr_{n,0} < \Ov{\vr}\ \mbox{uniformly in}\ n, \ \vr_{n,0} \to \vr_E(0, \cdot) \ \mbox{in}\ L^1 (\Omega), \br
		0 &< \underline{\vt} < \inf_{(0,T) \times \Omega} \vt_{n,0} \leq  \sup_{(0,T) \times \Omega} \vt_{n,0} < \Ov{\vt}\ \mbox{uniformly in}\ n,\ 
		e(\vr_{n,0}, \vt_{n,0})(0, \cdot) \to e_E(0, \cdot) \ \mbox{in}\ L^1(\Omega), \br
		&|\vu_{n,0} | \leq \Ov{\vu} \ \mbox{uniformly in}\ n,\ \vu_{n,0} \to \vu_E(0, \cdot) \ \mbox{in}\ L^1(\Omega; R^3).
		\nonumber
	\end{align}
	Finally, suppose there is a sequence $\delta_n \to 0$ such that 
	\begin{align}
		\frac{\mu_n}{\delta_n} &\to 0 \ \mbox{as}\ n \to \infty, \br
		\frac{1}{\delta_n} &\int_0^T \int_{\Omega_{\delta_n}} \vt^{2}_n \dx \dt \leq c,\br 
		\frac{1}{\mu_n}\int_0^T & \left\| \vr_n (\vu_n)_n \right\|_{L^{2}(\Omega_{\delta_n})}^2  \dt \to 0
		\label{r2bis}
	\end{align}
	as $n \to \infty$.	
	
	Then 
	\[
		\vr_n \to \vr_E,\ \vr_n e(\vr_n, \vt_n) \to \vr_E e_E  \ \mbox{in}\ L^1((0,T) \times \Omega),\ \vr_n \vu_n \to \vr_E \vu_E \ \mbox{in}\ L^1((0,T) \times \Omega; R^3). 
	\]
	
\end{Theorem}

The rest of the paper is devoted to the proof of the above results.

\section{Consistency of the vanishing dissipation/radiation approximation}
\label{c} 

As a preliminary step, we show \emph{consistency} of the vanishing dissipation/radiation approximation. 

\subsection{Temperature for the Euler system} \label{subsec:temp} 

First we introduce the temperature $\vt_E$ associated to the limit system. Without loss of generality, we may fix the constants $\underline{\vr}$, $\Ov{\vr}$ in \eqref{r2} so that 
\begin{equation} \label{c1} 
	0 < \underline{\vr} < \inf_{[0,T] \times \Ov{\Omega}} \vr_E \leq \sup_{[0,T] \times \Ov{\Omega}} \vr_E < \Ov{\vr}.
	\end{equation}
Next, in accordance with the hypotheses of Theorem \ref{Tr1}, 
\begin{equation} \label{c2} 
	0 < \underline{e} < \inf_{[0,T] \times \Ov{\Omega}} e_E \leq \sup_{[0,T] \times \Ov{\Omega}} e_E < \Ov{e}
\end{equation} 
for certain constants $\underline{e}$, $\Ov{e}$. Let us set
\begin{equation*} 
	\underline{Z} > \frac{\overline{\vr}}{\left((\gamma-1)\underline{e}\right)^{\frac{1}{\gamma-1}}} 
\end{equation*} 
and apply lemma \ref{lemma:P} to obtain suitable functions $P,S$. Furthermore we define 
\begin{equation*} 
	\vt_E := (\gamma-1) e_E.
\end{equation*}
Note that $\vt_E>(\gamma-1)\underline{e}$ and hence
$$
\frac{\vr_E}{(\vt_E)^{\frac{1}{\gamma-1}}}<\underline{Z}. 
$$ 
By virtue of \eqref{eq:lemmaP}, we have
\[
e(\vr_E, \vt_E) = e_E \ \mbox{in}\ [0,T] \times \Ov{\Omega}.
\] 

Moreover, without loss of generality, we may suppose 
\begin{equation} \label{c4} 
	0 < \underline{\vt} < \inf_{[0,T] \times \Ov{\Omega}} \vt_E \leq \sup_{[0,T] \times \Ov{\Omega}} \vt_E < \Ov{\vt},
\end{equation}
with the same constants $\underline{\vt}$, $\Ov{\vt}$ as in \eqref{r2}. From this moment on, the pressure law is fixed. 
	 
As $p$, $e$, and $s$ comply with Gibbs' relation, the smooth solution of the Euler system conserves the entropy:
\begin{equation} \label{c5}
	\partial_t (\vr_E s(\vr_E, \vt_E)) + \Div (\vr_E s(\vr_E, \vt_E) \vu_E ) = 0,
	\end{equation}
where $s$ is given by \eqref{i6}. 

\subsection{Consistency}

\label{Con}

The Navier--Stokes--Fourier system \eqref{p2} may be viewed a singular perturbation of the Euler system with the extra ``error'' terms 
\begin{align} 
E^1_n &= p_R =\frac{a_n}{3} \vt^4_n,\br E^2_n &= \mu_n \mathbb{S} (\vt_n, \Grad \vu_n) = \mu_n \left( \tmu (\vt_n) \left( \Grad \vu_n + \Grad^t \vu_n - \frac{2}{3} \Div \vu_n \mathbb{I} \right) 
+ \teta(\vt_n) \Div \vu_n \mathbb{I} \right), \br
E^3_n &= \vr s_R = \frac{4 a_n}{3} \vt_n^3,\ E^4_n = \vr s_R \vu = \frac{4 a_n}{3} \vt_n^3 \vu_n,\ E^5_n \equiv \kappa_n \frac{\vc{q}}{\vt} = \kappa_n \tkap (\vt_n) \frac{\Grad \vt_n}{\vt_n},
\
E^6_n = \vr e_R = a_n \vt_n^4.
\label{c6}
\end{align}

We say that the approximation of the Euler system by the Navier--Stokes--Fourier system is consistent, if the above ``error'' terms vanish in the asymptotic limit $n \to 0$.  
As a matter of fact, we need a milder form of consistency compatible with the relative energy inequality. More specifically, it is sufficient to control the ``errors'' by the dissipation term
\begin{align}
\mathcal{D}_n &\equiv \mu_n \intO{ \frac{ \tmu (\vt_n) }{\vt_n } \left| \Grad \vu_n + \Grad^t \vu_n - \frac{2}{3} \Div \vu_n \mathbb{I} \right|^2  }
+ \mu_n \intO{ \frac{ \teta (\vt_n) }{\vt_n } |\Div \vu|^2 } \br 
&+ \kappa_n \intO{ \frac{\tkap (\vt_n) }{\vt_n^2} |\Grad \vt_n |^2 },
\nonumber
\end{align}
and the total energy 
\[
\mathcal{E}_n \equiv \intO{ \left( \frac{1}{2} \vr_n |\vu_n|^2 + \vr_n e(\vr_n, \vt_n) + a_n \vt^4_n \right) }.
\]
For each error term $E^i_n$, $i=1,\dots, 6$ specified in \eqref{c6} and $\ep > 0$, we have to find $c(\ep)$ such that 
\begin{equation} \label{c7}
\| E^i_n \|_{L^1(\Omega)} \leq \ep \mathcal{D}_n + c(\ep) \mathcal{E}_n + c(\ep) \omega_n \ \mbox{uniformly for}\ n \to \infty,\ \omega_n \to 0. 
\end{equation}
Obviously, $E^1_n = p_R$, $E^6_n = \vr e_R$, and $E^3_n = \vr s_R$ satisfy \eqref{c7} (with $\ep = 0$), it remains to handle the viscous stress, the heat flux and the entropy convective flux term.

Moreover, we recall some basic estimates that follow directly from the hypotheses \eqref{i6}, \eqref{p5}, \eqref{eq:prop-SP1}--\eqref{eq:prop-SP4}: 
\begin{align} 
	\vr^{\gamma} + \vr \vt &\aleq \vr e(\vr, \vt), \br
	0 \leq \vr s(\vr, \vt) &\aleq \vr \left( 1 + |\log(\vr)| + [\log(\vt)]^+ \right).
	\label{c7a}
\end{align}

\subsubsection{Viscous stress consistency}

By virtue of hypothesis \eqref{p3}, 
\begin{align} 
	\int_\Omega &\mu_n \tmu (\vt_n) \left| \Grad \vu_n + \Grad^t \vu_n - \frac{2}{3} \Div \vu_n \mathbb{I} \right| \dx \br
	&\leq \ep \mu_n \intO{ \frac{ \tmu (\vt_n) }{\vt_n } \left| \Grad \vu_n + \Grad^t \vu_n - \frac{2}{3} \Div \vu_n \mathbb{I} \right|^2  } + c(\ep) \mu_n \intO{ \left(1 + \vt_n^{1 + \alpha} \right) } \br
	&\leq \ep \mathcal{D}_n + c(\ep) \mu_n + c(\ep) \mu_n \intO{ \left[\vt_n^{1 + \alpha} \right]_{\rm res} } \br
	&\leq \ep \mathcal{D}_n + c(\ep) \mu_n + c(\ep) \frac{\mu_n}{a_n^{\frac{1+\alpha}{4}}} \intO{ \left[a_n^{\frac{1+\alpha}{4}} \vt_n^{1 + \alpha} \right]_{\rm res} } \br
	&\leq \ep \mathcal{D}_n + c(\ep) \mu_n + c(\ep) \intO{ \left[a_n \vt_n^4 +1 \right]_{\rm res} }
	\nonumber
\end{align}
where the last inequality follows from hypothesis \eqref{r1} and the simple fact that $x^{1+\alpha}\leq c (x^4 + 1)$. Thus we obtain the desired estimate \eqref{c7}. The bulk viscosity term can be handled in a similar fashion. 

\subsubsection{Heat flux consistency}

Similarly to the preceding part, 
\begin{align}
	\int_\Omega &\kappa_n \frac{\tkap (\vt_n) }{\vt_n} |\Grad \vt_n | \dx \leq \ep \kappa_n \intO{ \frac{\tkap (\vt_n) }{\vt_n^2} |\Grad \vt_n |^2 } + 
	c(\ep) \kappa_n \intO{ \tkap(\vt_n) } \br 
	&\leq \ep \mathcal{D}_n + c(\ep) \kappa_n + c(\ep) \kappa_n \intO{ \vt_n^3 } \leq \ep \mathcal{D}_n + c(\ep) \kappa_n + c(\ep) \frac{\kappa_n}{a_n^{\frac{3}{4}}} \left( \intO{ a_n \vt_n^4 } \right)^{\frac{3}{4}} \br 
	&\leq \ep \mathcal{D}_n + c(\ep) \kappa_n + c(\ep) \frac{\kappa_n}{a_n^{\frac{3}{4}}} 
	+ c(\ep) \mathcal{E}_n; 
	\label{c8b}
\end{align}
whence \eqref{c7} follows from hypothesis \eqref{r1}. 

\subsubsection{Radiation entropy convective flux consistency}

To close the circle of consistence estimates, we have to handle the integral
\[
a_n \intO{ \vt^3_n |\vu_n | } \leq a_n \| \vt_n^3 \|_{L^{\frac{4}{3}}(\Omega)} \| \vu_n \|_{L^4(\Omega; R^3)}
\]
that corresponds to the radiation entropy convective flux.

By virtue of Sobolev embedding theorem, 
\[
\| \vu_n \|_{L^4(\Omega; R^3)} \aleq \| \vu_n \|_{W^{1, \frac{8}{5 - \alpha}}(\Omega; R^3)} \ \mbox{as long as}\ \alpha \geq \frac{1}{3},
\]
and, by a generalized Korn-Poincar\' e inequality \cite[Theorem 11.23]{FeNo6A}, 
\[
\| \vu_n \|_{W^{1, \frac{8}{5 - \alpha}}(\Omega; R^3)} \aleq \left(  \left\| \Grad \vu_n + \Grad^t \vu_n - \frac{2}{3} \Div \vu_n \mathbb{I} 
\right\|_{L^{\frac{8}{5 - \alpha}}(\Omega; R^9)}   + \intO{\vr_n |\vu_n| }                  \right).
\]
Another application of H\" older's inequality yields
\begin{align}
\left\| \Grad \vu_n + \Grad^t \vu_n - \frac{2}{3} \Div \vu_n \mathbb{I} 
\right\|_{L^{\frac{8}{5 - \alpha}}(\Omega; R^9)} \br \leq \| \vt_n^{\frac{1 - \alpha}{2}} \|_{L^{\frac{8}{1 - \alpha}}(\Omega)} 
\left\| \vt_n^{\frac{\alpha - 1}{2}} \left( \Grad \vu_n + \Grad^t \vu_n - \frac{2}{3} \Div \vu_n \mathbb{I} \right)
\right\|_{L^{2}(\Omega; R^9)} 
.
\nonumber
\end{align}
Consequently, 
\begin{align}
a_n &\intO{ \vt^3_n |\vu_n | } \leq a_n \| \vt_n^3 \|_{L^{\frac{4}{3}}(\Omega)} \| \vu_n \|_{L^4(\Omega; R^3)} \aleq 
 a_n \| \vt_n^3 \|_{L^{\frac{4}{3}}(\Omega)} \intO{\vr_n |\vu_n| } \br	
&+ a_n \| \vt_n^3 \|_{L^{\frac{4}{3}}(\Omega)}\| \vt_n^{\frac{1 - \alpha}{2}} \|_{L^{\frac{8}{1 - \alpha}}(\Omega)} 
\left\| \vt_n^{\frac{\alpha - 1}{2}} \left( \Grad \vu_n + \Grad^t \vu_n - \frac{2}{3} \Div \vu_n \mathbb{I} \right)
\right\|_{L^{2}(\Omega; R^9)} \br 
&\leq \ep \mathcal{D}_n + c(\ep) \frac{a_n^2}{\mu_n}  \| \vt_n^3 \|_{L^{\frac{4}{3}}(\Omega)}^2 \| \vt_n^{\frac{1 - \alpha}{2}} \|_{L^{\frac{8}{1 - \alpha}}(\Omega)}^2 + 
c a_n \| \vt_n^3 \|_{L^{\frac{4}{3}}(\Omega)} \intO{\vr_n |\vu_n| }.
\label{c8}
	\end{align}

Finally, by virtue of hypothesis \eqref{r1}
\[
\frac{a_n^2}{\mu_n}  \| \vt_n^3 \|_{L^{\frac{4}{3}}(\Omega)}^2 \| \vt_n^{\frac{1 - \alpha}{2}} \|_{L^{\frac{8}{1 - \alpha}}(\Omega)}^2 = 
\frac{a_n^2}{\mu_n} \left( \intO{ \vt_n^4 } \right)^{\frac{7 - \alpha}{4}} = 
\frac{a_n^{\frac{1 + \alpha}{4}} }{\mu_n} \left( \intO{ a_n \vt_n^4 } \right)^{\frac{7 - \alpha}{4}} 
\aleq \mathcal{E}_n^{\frac{ 7 - \alpha }{4}}.
\]
The rightmost integral in \eqref{c8} can be handled in a similar fashion. Since $\alpha\in [\frac{1}{3},1]$, we have in particular $0 \leq \alpha \leq 3$ and the desired conclusion \eqref{c7} follows from boundedness of the total energy.

\section{Convergence}
\label{C}

The proof of convergence consists in plugging the strong solution $(\vr_E, \vt_E, \vu_E)$ of the Euler system as the test functions $r = \vr_E$, $\Theta = \vt_E$, $
\vU = \vu_E$ in the relative energy inequality \eqref{p11}. This can be done in a direct manner in the case of the complete slip boundary conditions \eqref{i13} , whereas the velocity $\vu_E$ must be modified to comply with 
the homogeneous Dirichlet boundary conditions in the case of no--slip \eqref{i14}. We focus on the latter case as the proof in the case of the complete slip boundary conditions can be performed in a way similar to \cite{Fei2015A}.

\subsection{Velocity regularization} \label{subsection:C1}

If the solutions of the Navier--Stokes--Fourier system satisfy the no--slip boundary conditions, the velocity $\vu_E$ is not eligible for the relative energy inequality \eqref{p11} as its tangential component may not vanish on $\partial \Omega$. Instead we consider
\begin{equation} \label{C1}
\vU = \vu_E - \vc{v}_\delta, 
\end{equation} 
where the perturbation $\vc{v}_\delta$ is given as 
\begin{equation} \label{C2}
\vc{v}_\delta (t,x) = \xi \left( \frac{ {\rm dist}[x, \partial \Omega] }{\delta} \right) \vu_E (t, \Pi(x)), \ \delta > 0,
\end{equation}	
where 
\[ 
\xi \in C^\infty(R),\ \xi' \leq 0, \ \xi (d) = 1 \ \mbox{if}\ d \leq 0, \ \xi(d) = 0 \ \mbox{if}\ \xi \geq 1,
\]
and 
\[
\Pi(x) \in \partial \Omega \ \mbox{is the nearest point to} \ x \ \mbox{in} \ \partial \Omega.
\]

If $\partial \Omega$ is of class $C^k$, $k \geq 2$, then ${\rm dist}[x, \partial \Omega] \in C^k(\Omega_\delta)$ for any $0 < \delta < \delta_0$, and
\[
\Grad {\rm dist}[x, \partial \Omega] = \frac{ x - \Pi(x) }{|x - \Pi(x)|} = - \vc{n} (\Pi(x)) \ \mbox{for any}\ x \in \Omega_\delta,
\]
see Foote \cite{FO}.

\subsection{Application of the relative energy inequality}

As $\vc{U} = \vu_E - \vd$ vanishes on $\partial \Omega$, the trio $(r = \vr_E, \vc{U} = \vu_E - \vd, \Theta = \vt_E)$ can be used as test functions in the relative energy inequality 
\eqref{p11}. Recall that at this stage we have the following vanishing parameters: $\mu_n$, $\kappa_n$, $a_n$, and $\delta = \delta_n$. 

We have 
\begin{equation}
\left| E \left(\vr_n, \vt_n, \vu_n \Big| \vr_E, \vt_E, \vu_E \right) - E \left(\vr, \vt, \vu \Big| \vr_E, \vt_E , \vu_E - \vd \right) \right| \aleq
\left| \vr_n (\vu_n - \vu_E) \vd \right| + \vr_n |\vd|^2
\label{C2a}
\end{equation}
Seeing that 
\begin{equation}
{\rm ess} \sup_{t \in (0,T) } \| \vr_n \|_{L^\gamma (\Omega)} + 
{\rm ess} \sup_{t \in (0,T) } \| \vr_n \vu_n \|_{L^{\frac{2 \gamma}{\gamma + 1}} (\Omega; R^3)} \aleq 1,
\label{C2b}
\end{equation}
we may infer that 
\begin{equation} \label{C2c}
\intO{ \left|  E \left(\vr_n, \vt_n, \vu_n \Big| \vr_E, \vt_E, \vu_E \right) - E \left(\vr, \vt, \vu \Big| \vr_E, \vt_E , \vu_E - \vd \right) \right|} \to 0 
\ \mbox{as}\ \delta \to 0.
\end{equation}

The first rather straightforward observation 
is that, under hypothesis \eqref{r2} concerning the initial data,
\[
\intO{ {E}_{a_n} \left(\vr_{0,n}, \vt_{0,n}, \vu_{0,n} \Big| \vr_E(0, \cdot), \vt_E (0, \cdot), \vu_E(0, \cdot) - \vd(0, \cdot) \right) } \to 0 
\ \mbox{for}\ n \to \infty,\ \delta \to 0.
\]	
Consequently, we can write \eqref{p11} in the form
\begin{align}
	&\intO{ E_{a_n} \left(\vr_n,\vt_n,\vu_n \Big| \vr_E, \vt_E, \vu_E - \vd \right)(\tau, \cdot) } \br &+ {\underline{\vt}} \int_0^\tau \int_\Omega \frac{1}{\vt_n}\left( \mu_n \S(\vt_n,\Grad \vu_n):\Grad\vu_n - \kappa_n \frac{\vq(\vt_n,\Grad\vt_n)\cdot \Grad\vt_n}{\vt_n}\right)\dx\dt \br
	&\leq\int_0^\tau \int_\Omega \vr_n(\vu_n - (\vu_E - \vd))\cdot\Grad (\vu_E - \vd) \cdot (\vu_n - (\vu_E - \vd)) \dx\dt  \br &\quad + \mu_n \int_0^\tau \int_\Omega \S(\vt_n,\Grad\vu_n) :\Grad 
	(\vu_E - \vd) \dx\dt 
	- \kappa_n \int_0^\tau \int_\Omega \frac{\vq(\vt_n,\Grad\vt_n)}{\vt_n}\cdot \Grad\vt_E \dx\dt \br &\quad+ \int_0^\tau \int_\Omega \vr_n \big((s + s_R) (\vr_n,\vt_n) - (s + s_R)(\vr_E,\vt_E)\big)((\vu_E - \vd)-\vu_n)\cdot \Grad\vt_E \dx\dt \br
	&\quad +\int_0^\tau\int_\Omega \vr_n \Big(\partial_t (\vu_E - \vd) + (\vu_E - \vd) \cdot \Grad (\vu_E - \vd) \Big)\cdot((\vu_E - \vd)-\vu_n)\dx\dt \br &\quad - \int_0^\tau\int_\Omega (p + p_R)(\vr_n,\vt_n)\Div (\vu_E - \vd) \dx\dt \br
	&\quad - \int_0^\tau\int_\Omega \vr_n \big((s + s_R) (\vr_n,\vt_n )-(s + s_R) (\vr_E,\vt_E)\big)\Big(\partial_t\vt_E + (\vu_E - \vd) \cdot \Grad\vt_E \Big)\dx\dt \br 
	&\quad + \int_0^\tau\int_\Omega \left(\left(1-\frac{\vr_n}{\vr_E}\right)\partial_t (p + p_R) (\vr_E,\vt_E) - \frac{\vr_n}{\vr_E}\vu_n \cdot \Grad (p + p_R) (\vr_E,\vt_E)\right)\dx\dt 
	+ h (n, \delta) ,\label{C3}
\end{align} 
holds for a.a. $\tau \in (0,T)$, where $h$ denotes a generic sequence, 
\[
h(n, \delta) \to 0 \ \mbox{as}\ n \to \infty,\ \delta \to 0.
\]

Our goal is to show 
\[
\intO{ E_{a_n} \left(\vr_n,\vt_n,\vu_n \Big| \vr_E, \vt_E, \vu_E - \vd \right)(\tau, \cdot) } = h(n, \delta) \ \mbox{uniformly for a.a.}\ \tau \in (0,T),
\]
by means of a Gronwall type argument.

\subsection{Integrals controlled by the consistency estimates}

Evoking the bounds obtained in Section \ref{Con} we get 
\begin{align}
&\kappa_n \left| \int_\Omega \frac{\vq(\vt_n,\Grad\vt_n)}{\vt_n}\cdot \Grad\vt_E \dx \right| = 
\kappa_n \left| \intO{ \frac{ \tkap{(\vt_n) } } {\vt_n}  \Grad \vt_n \cdot \Grad \vt_E } \right| \br 
&\quad \leq \kappa_n \left| \intO{ \left[ \frac{ \tkap{(\vt_n) } } {\vt_n} \right]_{\rm ess}  \Grad \vt_n \cdot \Grad \vt_E } \right| + 
\kappa_n \left| \intO{ \left[ \frac{ \tkap{(\vt_n) } } {\vt_n} \right]_{\rm res}  \Grad \vt_n \cdot \Grad \vt_E } \right| \br 
&\quad \leq \ep \mathcal{D}_n + c(\ep, \| \Grad \vt_E \|_{L^\infty} ) \kappa_n + \| \Grad \vt_E \|_{L^\infty} 
\kappa_n  \intO{ \left[ \frac{ \tkap{(\vt_n) } } {\vt_n} \right]_{\rm res}  |\Grad \vt_n| }, 
\nonumber
\end{align}
where, by virtue of \eqref{c8b},
\begin{align}
\kappa_n  \intO{ \left[ \frac{ \tkap{(\vt_n) } } {\vt_n} \right]_{\rm res}  |\Grad \vt_n| } &\leq 
\ep \mathcal{D}_n + c(\ep) \intO{ a_n [ \vt_n^4 ]_{\rm res} } \br &\leq 
\ep \mathcal{D}_n + c(\ep) \intO{E_{a_n} \left(\vr_n,\vt_n,\vu_n \Big| \vr_E, \vt_E, \vu_E - \vd \right)}. 
\nonumber
\end{align}

Using the consistency estimates of Section \ref{Con}, we can handle other integrals containing vanishing parameters. Accordingly, the inequality \eqref{C3} simplifies to
\begin{align}
	&\intO{ E_{a_n} \left(\vr_n,\vt_n,\vu_n \Big| \vr_E, \vt_E, \vu_E - \vd \right)(\tau, \cdot) } \br &+ {\underline{\vt}} \int_0^\tau \int_\Omega \frac{1}{\vt_n}\left( \mu_n \S(\vt_n,\Grad \vu_n):\Grad\vu_n - \kappa_n \frac{\vq(\vt_n,\Grad\vt_n)\cdot \Grad\vt_n}{\vt_n}\right)\dx\dt \br
	&\leq - \int_0^\tau \int_\Omega \vr_n(\vu_n - (\vu_E - \vd))\cdot\Grad \vd \cdot (\vu_n - (\vu_E - \vd)) \dx\dt  \br &\quad - \mu_n \int_0^\tau \int_\Omega \S(\vt_n,\Grad\vu_n) :\Grad 
	\vd \dx\dt 
 \br &\quad+ \int_0^\tau \int_\Omega \vr_n \big((s + s_R) (\vr_n,\vt_n) - (s + s_R)(\vr_E,\vt_E)\big)((\vu_E - \vd)-\vu_n)\cdot \Grad\vt_E \dx\dt \br
	&\quad +\int_0^\tau\int_\Omega \vr_n \Big(\partial_t (\vu_E - \vd) + (\vu_E - \vd) \cdot \Grad (\vu_E - \vd) \Big)\cdot((\vu_E - \vd)-\vu_n)\dx\dt \br &\quad - \int_0^\tau\int_\Omega (p + p_R)(\vr_n,\vt_n)\Div (\vu_E - \vd) \dx\dt \br
	&\quad - \int_0^\tau\int_\Omega \vr_n \big(s  (\vr_n,\vt_n )- s  (\vr_E,\vt_E)\big)\Big(\partial_t\vt_E + (\vu_E - \vd) \cdot \Grad\vt_E \Big)\dx\dt \br 
	&\quad + \int_0^\tau\int_\Omega \left(\left(1-\frac{\vr_n}{\vr_E}\right)\partial_t p  (\vr_E,\vt_E) - \frac{\vr_n}{\vr_E}\vu_n \cdot \Grad p  (\vr_E,\vt_E)\right)\dx\dt \br
	&\quad + c \int_0^\tau \intO{ E_{a_n} \left(\vr_n,\vt_n,\vu_n \Big| \vr_E, \vt_E, \vu_E - \vd \right) } \dt   + h (n, \delta) ,\label{C4}
\end{align} 

Moreover, as $\vu_E \cdot \vc{n}|_{\partial \Omega} = 0$, 
\begin{equation} \label{C5}
	\|\Div \vd \|_{L^\infty} \aleq 1 \ \mbox{independently of}\ \delta,
\end{equation} 
and, consequently,
\begin{align}
	&\intO{ E_{a_n} \left(\vr_n,\vt_n,\vu_n \Big| \vr_E, \vt_E, \vu_E - \vd \right)(\tau, \cdot) } \br &+ {\underline{\vt}} \int_0^\tau \int_\Omega \frac{1}{\vt_n}\left( \mu_n \S(\vt_n,\Grad \vu_n):\Grad\vu_n - \kappa_n \frac{\vq(\vt_n,\Grad\vt_n)\cdot \Grad\vt_n}{\vt_n}\right)\dx\dt \br
	&\leq - \int_0^\tau \int_\Omega \vr_n \vu_n \cdot\Grad \vd \cdot (\vu_n - (\vu_E - \vd)) \dx\dt  \br &\quad - \mu_n \int_0^\tau \int_\Omega \S(\vt_n,\Grad\vu_n) :\Grad 
	\vd \dx\dt 
	\br &\quad+ \int_0^\tau \int_\Omega \vr_n [ s(\vr_n,\vt_n) + 1 ]_{\rm res}|\vu_n| |\Grad\vt_E| \dx\dt \br
	&\quad + \int_0^\tau \intO{ \vr_n \Big( \partial_t \vu_E + \vu_E \cdot \Grad \vu_E + \frac{1}{\vr_E} \Grad p(\vr_E, \vt_E) \Big) \cdot (\vu_E - \vu_n) } \dt \br
	&\quad - \int_0^\tau \intO{ \vr_n \Big( \partial_t \vd + \vd \cdot \Grad \vu_E \Big) \cdot ((\vu_E - \vd)-\vu_n) } \dt \br 
	 &\quad + \int_0^\tau\int_\Omega \Big( p(\vr_E, \vt_E) - p(\vr_n,\vt_n) \Big) \Div \vu_E  \dx\dt \br
	&\quad - \int_0^\tau\int_\Omega \vr_n \big(s  (\vr_n,\vt_n )- s  (\vr_E,\vt_E)\big)\Big(\partial_t\vt_E + \vu_E  \cdot \Grad\vt_E \Big)\dx\dt \br 
	&\quad + \int_0^\tau\int_\Omega \left(1-\frac{\vr_n}{\vr_E}\right) \left( \partial_t p  (\vr_E,\vt_E) + \vu_E \cdot \Grad p  (\vr_E,\vt_E)\right)\dx\dt \br
	&\quad + c \int_0^\tau \intO{ E_{a_n} \left(\vr_n,\vt_n,\vu_n \Big| \vr_E, \vt_E, \vu_E - \vd \right) } \dt  + h (n, \delta) .\label{C6}
\end{align} 

Finally, as $(\vr_E, \vt_E, \vu_E)$ solves the Euler system, 
\[
\partial_t \vu_E + \vu_E \cdot \Grad \vu_E + \frac{1}{\vr_E} \Grad p(\vr_E, \vt_E) = 0.
\]
In addition, it is easy to check that 
\begin{equation} \label{C7}
	\|\partial_t \vd \|_{L^\infty} +  \| \vd \|_{L^\infty} \aleq 1 \ \mbox{independently of}\ \delta.
	\end{equation} 
Consequently, 
\begin{align}
\int_0^\tau &\intO{ \vr_n \Big( \partial_t \vd + \vd \cdot \Grad \vu_E \Big) \cdot ((\vu_E - \vd)-\vu_n) } \dt \br &= 
\int_0^\tau \int_{\Omega_\delta} \vr_n \Big( \partial_t \vd + \vd \cdot \Grad \vu_E \Big) \cdot ((\vu_E - \vd)-\vu_n) \dx \dt \to 0 \ \mbox{as}\ \delta \to 0
\nonumber
\end{align}
as both $(\vr_n )_{n \geq 0}$ and $(\vr_n \vu_n)_{n \geq 0}$ are equi--integrable in $(0,T) \times \Omega$.
Thus \eqref{C6} reduces to 
\begin{align}
	&\intO{ E_{a_n} \left(\vr_n,\vt_n,\vu_n \Big| \vr_E, \vt_E, \vu_E - \vd \right)(\tau, \cdot) } \br &+ {\underline{\vt}} \int_0^\tau \int_\Omega \frac{1}{\vt_n}\left( \mu_n \S(\vt_n,\Grad \vu_n):\Grad\vu_n - \kappa_n \frac{\vq(\vt_n,\Grad\vt_n)\cdot \Grad\vt_n}{\vt_n}\right)\dx\dt \br
	&\leq - \int_0^\tau \int_\Omega \vr_n \vu_n \cdot\Grad \vd \cdot (\vu_n - (\vu_E - \vd)) \dx\dt  \br &\quad - \mu_n \int_0^\tau \int_\Omega \S(\vt_n,\Grad\vu_n) :\Grad 
	\vd \dx\dt 
	\br &\quad+ \int_0^\tau \int_\Omega \vr_n [ s(\vr_n,\vt_n) + 1 ]_{\rm res}|\vu_n| |\Grad\vt_E| \dx\dt \br 
	&\quad + \int_0^\tau\int_\Omega \Big( p(\vr_E, \vt_E) - p(\vr_n,\vt_n) \Big) \Div \vu_E  \dx\dt \br
	&\quad - \int_0^\tau\int_\Omega \vr_n \big(s  (\vr_n,\vt_n )- s  (\vr_E,\vt_E)\big)\Big(\partial_t\vt_E + \vu_E  \cdot \Grad\vt_E \Big)\dx\dt \br 
	&\quad + \int_0^\tau\int_\Omega \left(1-\frac{\vr_n}{\vr_E}\right) \left( \partial_t p  (\vr_E,\vt_E) + \vu_E \cdot \Grad p  (\vr_E,\vt_E)\right)\dx\dt \br
	&\quad + c \int_0^\tau \intO{ E_{a_n} \left(\vr_n,\vt_n,\vu_n \Big| \vr_E, \vt_E, \vu_E - \vd \right) }\dt   + h (n, \delta) .\label{C8}
\end{align} 

\subsection{Integrals independent of the boundary layer}

Now, we estimate the integrals on the right--hand side of \eqref{C8} that are independent of $\vd$. First, by virtue of \eqref{c7a}, 
\begin{align}
\left| \intO{ \vr_n [ s(\vr_n,\vt_n) + 1 ]_{\rm res}|\vu_n| |\Grad\vt_E| } \right|
&\aleq \intO{ [\vr_n]_{\rm res} |\vu_n|^2 } + \intO{ [\vr_n]_{\rm res} s^2(\vr_n, \vt_n) } \br 
&\aleq \intO{ [ E_{a_n} (\vr_n, \vt_n, \vu_n ) ]_{\rm res} } \br
&\aleq \intO{ E_{a_n} \left(\vr_n,\vt_n,\vu_n \Big| \vr_E, \vt_E, \vu_E - \vd \right) } + h(\delta).	
\label{C12}
\end{align}
We point out that this step depends in an essential way on the fact that $s$ satisfies the Third law of thermodynamics.

Next, we recall two identities that follow from the specific form of EOS \eqref{i3}, \eqref{i4}, namely 
\begin{align}
	\partial_t \vt_E + \vu_E \cdot \Grad \vt_E &= - (\gamma - 1) \vt_E \Div \vu_E, \br
	\partial_t p(\vr_E, \vt_E) + \vu_E \cdot \Grad p(\vr_E, \vt_E) &= - \gamma p(\vr_E, \vt_E) \Div \vu_E.
\nonumber
	\end{align}
Consequently, we get 
\begin{align}
\int_\Omega &\Big( p(\vr_E, \vt_E) - p(\vr_n,\vt_n) \Big) \Div \vu_E  \dx   - \int_\Omega \vr_n \big(s  (\vr_n,\vt_n )- s  (\vr_E,\vt_E)\big)\Big(\partial_t\vt_E + \vu_E  \cdot \Grad\vt_E \Big)\dx  \br 
&+\int_\Omega \left(1-\frac{\vr_n}{\vr_E}\right) \left( \partial_t p  (\vr_E,\vt_E) + \vu_E \cdot \Grad p  (\vr_E,\vt_E)\right)\dx	\br
&= \intO{ \Div\vu_E \Big[ p(\vr_E, \vt_E) - p(\vr_n,\vt_n) + (\gamma - 1) \vr_E \vt_E  \big(s  (\vr_n,\vt_n )- s  (\vr_E,\vt_E)\big) \Big] } \br 
&- \gamma \intO{ \Div \vu_E \left(1-\frac{\vr_n}{\vr_E}\right) p(\vrE, \vt_E)	} \br 
&+ (\gamma - 1) \intO{\vt_E (\vr_n - \vr_E) \big(s  (\vr_n,\vt_n )- s  (\vr_E,\vt_E)\big) \Div \vu_E } 
	\label{C13}
	\end{align}

Finally, we use the identity 
\begin{align} 
\gamma + \left( \frac{\vr_n}{\vr_E} - 1 \right) \gamma p(\vr_e, \vt_E) + 
\left( \frac{\partial p(\vr_E, \vt_E) }{\partial \vr} (\vr_E - \vr_n ) + \frac{\partial p(\vr_E, \vt_E) }{\partial \vt} (\vt_E - \vt_n ) \right) \br
+ (\gamma - 1) \vr_E \vt_E	\left( \frac{\partial s(\vr_E, \vt_E) }{\partial \vr} (\vr_E - \vr_n ) + \frac{\partial s(\vr_E, \vt_E) }{\partial \vt} (\vt_E - \vt_n ) \right) .
	\label{C14}
	\end{align}
Plugging \eqref{C14} to \eqref{C13} yields the desired estimate. Thus \eqref{C8} reduces to 
\begin{align}
	&\intO{ E_{a_n} \left(\vr_n,\vt_n,\vu_n \Big| \vr_E, \vt_E, \vu_E - \vd \right)(\tau, \cdot) } \br &+ {\underline{\vt}} \int_0^\tau \int_\Omega \frac{1}{\vt_n}\left( \mu_n \S(\vt_n,\Grad \vu_n):\Grad\vu_n - \kappa_n \frac{\vq(\vt_n,\Grad\vt_n)\cdot \Grad\vt_n}{\vt_n}\right)\dx\dt \br
	&\leq - \int_0^\tau \int_\Omega \vr_n \vu_n \cdot\Grad \vd \cdot (\vu_n - (\vu_E - \vd)) \dx\dt  \br &\quad - \mu_n \int_0^\tau \int_\Omega \S(\vt_n,\Grad\vu_n) :\Grad 
	\vd \dx\dt\br &\quad  + c \int_0^\tau \intO{ E_{a_n} \left(\vr_n,\vt_n,\vu_n \Big| \vr_E, \vt_E, \vu_E - \vd \right) }\dt   + h (n, \delta) .\label{C15}
\end{align}
Note that inequality \eqref{C15} already completes the proof of Theorem \ref{Tr1}, where we may take $\vd = 0$. 

\subsection{Boundary layer}

It remains to control the first two integrals on the right--hand side of \eqref{C15} that represent the effect of the boundary layer. 

\subsubsection{Viscous stress}

Similarly to Section \ref{Con}, we have 
\[ 
\mu_n	\left| \intO{ \mathbb{S} (\vt_n, \Grad \vu_n ) : \Grad \vd } \right| \leq \ep \mathcal{D}_n + c(\ep) \mu_n \intO{ \vt_n \left( 1 + \vt_n^\alpha \right) |\Grad \vd |^2 }, 
	\]
where 
\[
\mu_n \intO{ \vt_n \left( 1 + \vt_n^\alpha \right) |\Grad \vd |^2 } \aleq \frac{\mu_n}{\delta^2} \int_{\Omega_\delta} \left(1 + \vt^{1 + \alpha}_n \right) \ \dx 
\aleq \frac{\mu_n} {\delta} \left( 1 +  \frac{1}{\delta} \int_{\Omega_\delta} \vt^{1 + \alpha}_n \ \dx \right).
\]
Consequently, hypothesis \eqref{r2b} yields the desired estimate. Note that this type of estimates forces us to consider the thinckeness $\delta$ of the boundary layer asymptotically large than 
$\mu$, 
\[
\frac{\mu_n}{\delta_n} \to 0.
\]

Alternatively, following Sueur \cite{Sue1}, we may suppose \eqref{rr2b}, meaning 
\[
\sqrt{\mu_n} \| \mathbb{S} (\vt_n, \Grad \vu_n ) \|_{L^2((0,T) \times \Omega_{\mu_n}; R^d)} \to 0.
\]
Setting $\mu_n \approx \delta_n$, we get 
\begin{align}
\mu_n	& \left| \int_0^T \intO{ \mathbb{S} (\vt_n, \Grad \vu_n ) : \Grad \vd } \dt \right| \br &\leq 
\sqrt{\mu_n} \| \mathbb{S} (\vt_n, \Grad \vu_n ) \|_{L^2((0,T) \times \Omega_{\mu_n}; R^d)} \| \sqrt{\mu_n} \Grad \vc{v}_{\mu_n} \||_{L^2((0,T) \times \Omega_{\mu_n}; R^d)} \to 0.
\end{align}

\subsubsection{Convective term}

Finally, we consider 
\[
\int_\Omega \vr_n \vu_n \cdot\Grad \vd \cdot (\vu_n - (\vu_E - \vd)) \dx = 
\int_{\Omega_\delta} \vr_n \vu_n \cdot\Grad \vd \cdot (\vu_n - (\vu_E - \vd)) \dx \ \mbox{in}\ \Omega_\delta.
\]
Recall that 
\begin{align}
\vc{w}(t,x) &= \vc{w}_n(t,x) + \vc{w}_\tau(t,x), \br \vc{w}_n(t,x) &= \left( \vc{w} \cdot \Grad {\rm dist}[x, \partial \Omega] \right) \Grad {\rm dist}[x, \partial \Omega],\ 
\vc{w}_\tau(t,x) = \vc{w}(t,x) - \vc{w}_n(t,x).
\nonumber
\end{align}
Similarly, for a scalar function $F$, we decompose 
\[
\Grad F = \nabla_n F + \nabla_\tau F,\ \nabla_n F = \left( \Grad {\rm dist}[x, \partial \Omega] \cdot \Grad F \right)  \Grad {\rm dist}[x, \partial \Omega].
\]
In accordance with the definition of $\vd$, we get
\begin{equation}
	(\vd)_n = 0,\ \| \nabla_{\tau} \vd \|_{L^\infty} \aleq 1,\ \| \nabla_n \vd \|_{L^\infty} \aleq \frac{1}{\delta}.
	\label{C16}   
	\end{equation}

Now, 
\begin{align}
\int_{\Omega_\delta} &\vr_n \vu_n \cdot\Grad \vd \cdot (\vu_n - (\vu_E - \vd)) \dx \br &= \int_{\Omega_\delta} \vr_n (\vu_n)_n \cdot \Grad \vd \cdot (\vu_n - (\vu_E - \vd)) \dx + 
\int_{\Omega_\delta} \vr_n (\vu_n)_\tau \cdot\Grad \vd \cdot (\vu_n - (\vu_E - \vd)) \dx \br 	
&= \int_{\Omega_\delta} \vr_n (\vu_n)_n \cdot \nabla_n (\vd)_\tau \cdot (\vu_n - (\vu_E - \vd)) \dx + 
\int_{\Omega_\delta} \vr_n (\vu_n)_\tau \cdot\nabla_{\tau} \vd \cdot (\vu_n - (\vu_E - \vd)) \dx, 
\nonumber
\end{align}
where, by virtue of \eqref{C16}, 
\begin{align}
&\left| \int_{\Omega_\delta} \vr_n (\vu_n)_\tau \cdot\nabla_{\tau} \vd \cdot (\vu_n - (\vu_E - \vd)) \dx \right| \br &\quad \aleq 
\intO{ E \left(\vr_n,\vt_n,\vu_n \Big| \vr_E, \vt_E, \vu_E - \vd \right) } + \int_{\Omega_\delta} \vr_n |\vu_E - \vd || \vu_n - (\vu_E - \vd) | \dx.  
\nonumber
\end{align}
In view of \eqref{C2b}, $(\vr_n)_{n\geq 0}$, $(\vr_n \vu_n )_{n \geq 0}$ are equi--integrable; whence 
\[
\int_0^T \int_{\Omega_\delta} \vr_n |\vu_E - \vd || \vu_n - (\vu_E - \vd) | \dx \dt \to 0 \ \mbox{as}\ \delta \to 0
\]
uniformly in $n$.

Thus it remains to handle the integral
\[
\int_{\Omega_\delta} \vr_n (\vu_n)_n \cdot \nabla_n (\vd)_\tau \cdot (\vu_n - (\vu_E - \vd)) \dx.
\]
By H\" older's inequality and \eqref{C16}, 
\begin{align}
&\left| \int_{\Omega_\delta} \vr_n (\vu_n)_n \cdot \nabla_n (\vd)_\tau \cdot (\vu_n - (\vu_E - \vd)) \dx \right| \br &\quad \leq \frac{1}{\delta} \left\| \vr_n (\vu_n)_n \right\|_{L^{\frac{24}{17 + 3\alpha}}(\Omega_\delta)}
\left\| \vu_n - (\vu_E - \vd) \right\|_{L^{\frac{24}{7 - 3\alpha}}(\Omega_\delta; R^3)},
\label{C17}
\end{align}
where $\frac{24}{7 - 3\alpha}$ is the critical exponent in the Sobolev--Poincar\' e inequality
\begin{equation} \label{C17bis}
\left\| \vu_n \right\|_{L^{\frac{24}{7 - 3\alpha}}(\Omega_\delta; R^3)} \aleq \| \Grad \vu_n \|_{L^{\frac{8}{5 - \alpha}}(\Omega_\delta; R^9)}.
\end{equation}
As $\vu_n|_{\partial \Omega} = 0$, Korn's inequality yields
\begin{align}
\left\| \vu_n \right\|_{L^{\frac{24}{7 - 3\alpha}}(\Omega_\delta; R^3)} &\aleq \| \Grad \vu_n \|_{L^{\frac{8}{5 - \alpha}}(\Omega_\delta; R^9)} \aleq 
\left\| \Grad \vu_n + \Grad \vu_n^t - \frac{2}{3} \Div \vu_n \mathbb{I} \right\|_{L^{\frac{8}{5 - \alpha}}(\Omega_\delta; R^9)} \br 
&\aleq \left\| \vt_n^{\frac{1 - \alpha}{2}} \right\|_{L^{\frac{8}{1 - \alpha}}(\Omega_\delta)} \left\| \vt^{\frac{\alpha - 1}{2}} \left( \Grad \vu_n + \Grad \vu_n^t - \frac{2}{3} \Div \vu_n \mathbb{I} \right) \right\|_{L^2(\Omega; R^9)}
\nonumber
\end{align}	 
Note that the constants are independent of $\delta$ as $\vu_n$ can be extended to be zero outside $\Omega$.

Thus going back to \eqref{C17} we deduce 
\begin{align}
&\left| \int_{\Omega_\delta} \vr_n (\vu_n)_n \cdot \nabla_n (\vd)_\tau \cdot (\vu_n - (\vu_E - \vd)) \dx \right| \br &\quad \leq 
\frac{c(\ep)}{\delta} \left\| \vr_n (\vu_n)_n \right\|_{L^{\frac{24}{17 + 3\alpha}}(\Omega_\delta)} 
+ \frac{c(\ep)}{\delta^2 \mu_n}	 \left\| \vr_n (\vu_n)_n \right\|_{L^{\frac{24}{17 + 3\alpha}}(\Omega_\delta)}^2 \left\| \vt_n^{\frac{1 - \alpha}{2}} \right\|_{L^{\frac{8}{1 - \alpha}}(\Omega_\delta)}^2 \br 
&\quad + \ep \mathcal{D}_n
\nonumber
	\end{align}
in accordance with hypothesis \eqref{r2b}. 

Finally, we consider $\alpha = 1$ and replace the critical exponent $\frac{24}{17 + 3 \alpha}$ in \eqref{C17}  by the $L^2$ norm. Consequently, 
\begin{align}
	&\left| \int_{\Omega_\delta} \vr_n (\vu_n)_n \cdot \nabla_n (\vd)_\tau \cdot (\vu_n - (\vu_E - \vd)) \dx \right| \br &\quad \leq \frac{1}{\delta} \left\| \vr_n (\vu_n)_n \right\|_{L^{2}(\Omega_\delta)}
	\left\| \vu_n - (\vu_E - \vd) \right\|_{L^{2}(\Omega_\delta; R^3)} \br 
	&\quad \aleq \frac{1}{\delta} \left\| \vr_n (\vu_n)_n \right\|_{L^{2}(\Omega_\delta)} \left\| \vu_n \right\|_{L^{2}(\Omega_\delta; R^3)} + 
	 \frac{1}{\sqrt{\delta}} \left\| \vr_n (\vu_n)_n \right\|_{L^{2}(\Omega_\delta)} \br 
	 &\quad \aleq \frac{1}{\delta} \left\| \vr_n (\vu_n)_n \right\|_{L^{2}(\Omega_\delta)} \left\| \vu_n \right\|_{L^{2}(\Omega_\delta; R^3)} + 
	 \sqrt{\frac{\mu_n}{\delta}} \left( 1 + \frac{1}{\mu_n} \left\| \vr_n (\vu_n)_n \right\|_{L^{2}(\Omega_\delta)}^2 \right).
	\label{C19}
\end{align}

Now, replacing \eqref{C17bis} by Hardy--Sobolev inequality, we gain multiplicative factor $\delta$, 
\begin{equation} \label{C17bos}
	\left\| \vu_n \right\|_{L^{2}(\Omega_\delta; R^3)} \aleq \delta \| \Grad \vu_n \|_{L^{2}(\Omega_\delta; R^9)}.
\end{equation}
Thus the final inequality reads 
\begin{align}
	&\left| \int_{\Omega_\delta} \vr_n (\vu_n)_n \cdot \nabla_n (\vd)_\tau \cdot (\vu_n - (\vu_E - \vd)) \dx \right| \br &\quad \leq c(\ep) 
	\left( \sqrt{\frac{\mu_n}{\delta}} \left( 1 + \frac{1}{\mu_n} \left\| \vr_n (\vu_n)_n \right\|_{L^{2}(\Omega_\delta)}^2 \right) + \frac{1}{\mu_n} \left\| \vr_n (\vu_n)_n \right\|_{L^{2}(\Omega_\delta)}^2
	\right) \br 
	&\quad + \ep \mathcal{D}_n
	\nonumber
\end{align}
in accordance with \eqref{r2bis}.

Alernatively, following Sueur \cite{Sue1}, we may suppose \eqref{rr2b}. The proof is exactly the same as in \cite{Sue1}.

\subsection{Strong convergence}
\label{sc}

We have established the convergence 
\[
\intO{ E_{a_n} \left( \vr_n, \vt_n, \vu_n \Big| \vr_E, \vt_E, \vu_E - \vc{v}_{\delta_n} \right) (\tau, \cdot) } \to 0 \ \mbox{as}\ n \to \infty 
\]
uniformly for a.a. $\tau \in (0,T)$. This obviously yields 
\[
\intO{ E \left( \vr_n, \vt_n, \vu_n \Big| \vr_E, \vt_E, \vu_E  \right) (\tau , \cdot) } \to 0 \ \mbox{as}\ n \to \infty. 
\]
In addition, as the energy of the initial data converges and both Euler and the Navier--Stokes--Fourier system conserve energy, we get 
\[
\intO{ \left( \frac{1}{2} \vr_n |\vu_n|^2 + \vr_n e(\vr_n, \vt_n) \right) } \to 
\intO{ \left( \frac{1}{2} \vr_E |\vu_E|^2 + \vr_E e(\vr_E, \vt_E) \right) } \ \mbox{in}\ L^1(0,T). 	
\]
This yields the desired strong convergence claimed in \eqref{r3}.

\def\cprime{$'$} \def\ocirc#1{\ifmmode\setbox0=\hbox{$#1$}\dimen0=\ht0
	\advance\dimen0 by1pt\rlap{\hbox to\wd0{\hss\raise\dimen0
			\hbox{\hskip.2em$\scriptscriptstyle\circ$}\hss}}#1\else {\accent"17 #1}\fi}


\begin{thebibliography}{10}
	
	\bibitem{BE}
	E.~Becker.
	\newblock {\em Gasdynamik}.
	\newblock Teubner-Verlag, Stuttgart, 1966.
	
	\bibitem{BEL1}
	F.~Belgiorno.
	\newblock Notes on the third law of thermodynamics, i.
	\newblock {\em J. Phys. A}, {\bf 36}:8165--8193, 2003.
	
	\bibitem{BEL2}
	F.~Belgiorno.
	\newblock Notes on the third law of thermodynamics, ii.
	\newblock {\em J. Phys. A}, {\bf 36}:8195--8221, 2003.
	
	\bibitem{BenSer}
	S.~Benzoni-Gavage and D.~Serre.
	\newblock {\em Multidimensional hyperbolic partial differential equations,
		{F}irst order systems and applications}.
	\newblock Oxford Mathematical Monographs. The Clarendon Press Oxford University
	Press, Oxford, 2007.
	
	\bibitem{ChauFei}
	N.~Chaudhuri and E.~Feireisl.
	\newblock {N}avier--{S}tokes--{F}ourier system with {D}irichlet boundary
	conditions.
	\newblock {\em {\bf arxiv preprint No. 2106.05315}}, 2021.
	

	
	\bibitem{Cowp}
	M.~Cowperthwaite.
	\newblock Relations between incomplete equations of state.
	\newblock {\em J. Franklin Institute}, {\bf 287}:379--387, 1969.
	
	\bibitem{Daf4}
	C.M. Dafermos.
	\newblock The second law of thermodynamics and stability.
	\newblock {\em Arch. Rational Mech. Anal.}, {\bf 70}:167--179, 1979.
	
	\bibitem{Fei2015A}
	E.~Feireisl.
	\newblock Vanishing dissipation limit for the {N}avier-{S}tokes-{F}ourier
	system.
	\newblock {\em Commun. Math. Sci.}, {\bf 14}(6):1535--1551, 2016.
	
	\bibitem{FeiNovOpen}
	E.~Feireisl and A.~Novotn{\' y}.
	\newblock {\em Mathematics of open fluid systems}.
	\newblock Birkh{\" a}user--Verlag, Basel.
	\newblock To appear.
	
	\bibitem{FeiNov10}
	E.~Feireisl and A.~Novotn{\' y}.
	\newblock Weak-strong uniqueness property for the full
	{N}avier-{S}tokes-{F}ourier system.
	\newblock {\em Arch. Rational Mech. Anal.}, {\bf 204}:683--706, 2012.
	
	\bibitem{FeNo6A}
	E.~Feireisl and A.~Novotn\'y.
	\newblock {\em Singular limits in thermodynamics of viscous fluids}.
	\newblock Advances in Mathematical Fluid Mechanics. Birkh\"auser/Springer,
	Cham, 2017.
	\newblock Second edition.
	
	\bibitem{FeiNov20}
	E.~Feireisl and A.~Novotn{\' y}.
	\newblock {N}avier-{S}tokes-{F}ourier system with general boundary conditions.
	\newblock {\em Arxive Preprint Series}, {\bf arXiv 2009.08207}, 2020.
	\newblock To appear in Commun. Math. Phys.
	
	\bibitem{FO}
	R.L. Foote.
	\newblock Regularity of the distance function.
	\newblock {\em Proc. Amer. Math. Soc.}, {\bf 92}:153--155, 1984.
	
	\bibitem{JesJiNov}
	D.~Jessl\'{e}, B.~J. Jin, and A.~Novotn\'{y}.
	\newblock Navier-{S}tokes-{F}ourier system on unbounded domains: weak
	solutions, relative entropies, weak-strong uniqueness.
	\newblock {\em SIAM J. Math. Anal.}, {\bf 45}(3):1907--1951, 2013.
	
	\bibitem{Kato}
	T.~Kato.
	\newblock Remarks on the zero viscosity limit for nonstationary
	{N}avier--{S}tokes flows with boundary.
	\newblock {\em In Seminar on PDE's, S.S. Chern (ed.), Springer, New York},
	1984.
	
	\bibitem{Kato1984}
	T.~Kato.
	\newblock Remarks on zero viscosity limit for nonstationary {N}avier-{S}tokes
	flows with boundary.
	\newblock In {\em Seminar on nonlinear partial differential equations
		({B}erkeley, {C}alif., 1983)}, volume~2 of {\em Math. Sci. Res. Inst. Publ.},
	pages 85--98. Springer, New York, 1984.
	
	\bibitem{MURU}
	I.~M{\" u}ller and T.~Ruggeri.
	\newblock {\em Rational extended thermodynamics}.
	\newblock Springer Tracts in Natural Philosophy 37, Springer-Verlag,
	Heidelberg, 1998.
	
	\bibitem{OX}
	J.~Oxenius.
	\newblock {\em Kinetic theory of particles and photons}.
	\newblock Springer-Verlag, Berlin, 1986.
	
	\bibitem{PokSkr}
	M.~Pokorn\'{y} and E.~Sk\v{r}\'{\i}\v{s}ovsk\'{y}.
	\newblock Weak solutions for compressible {N}avier-{S}tokes-{F}ourier system in
	two space dimensions with adiabatic exponent almost one.
	\newblock {\em Acta Appl. Math.}, 172:Paper No. 1, 31, 2021.
	
	\bibitem{SCHO1}
	S.~Schochet.
	\newblock The compressible {E}uler equations in a bounded domain: {E}xistence
	of solutions and the incompressible limit.
	\newblock {\em Commun. Math. Phys.}, {\bf 104}:49--75, 1986.
	
	\bibitem{Sue1}
	F.~Sueur.
	\newblock On the inviscid limit for the compressible {N}avier-{S}tokes system
	in an impermeable bounded domain.
	\newblock {\em J. Math. Fluid Mech.}, {\bf 16}(1):163--178, 2014.
	
	\bibitem{WangZhu}
	Y.-G. Wang and S.-Y. Zhu.
	\newblock On the vanishing dissipation limit for the full
	{N}avier-{S}tokes-{F}ourier system with non-slip condition.
	\newblock {\em J. Math. Fluid Mech.}, {\bf 20}(2):393--419, 2018.
	
\end{thebibliography}

\end{document}